\documentclass[11pt,a4paper,reqno]{amsart}
\usepackage{tikz}
\usepackage[utf8]{inputenc}
\usepackage{amsmath,bm,bbm,amsthm, amssymb}
\usepackage{fullpage}
\usepackage{color}
\usepackage{animate}
\usepackage{etoolbox}
\usepackage{comment}
\usepackage{pgf}
\usepackage{hyperref}
\usetikzlibrary{calc}
\usetikzlibrary{patterns}
\usetikzlibrary{arrows}
\usetikzlibrary{decorations.pathreplacing}
\usepackage[utf8]{inputenc}
\usepackage{dsfont}

\theoremstyle{plain}
\newtheorem{theorem}{Theorem}
\newtheorem{proposition}[theorem]{Proposition}
\newtheorem{corollary}[theorem]{Corollary}
\newtheorem{lemma}[theorem]{Lemma}

\theoremstyle{definition}
\newtheorem{definition}[theorem]{Definition}

\newtheorem{example}[theorem]{Example}

\theoremstyle{remark}

\usepackage{hyperref}

\definecolor{wiasblue}   {cmyk}{1.0, 0.60, 0, 0}
\definecolor{mlugreen}{RGB}{0,81,51}

\def\E{\mathbb E}
\def\P{\mathbb P}
\def\M{\mathbb M}

\def\R{\mathbb R}

\def\ms{\mathsf}

\newcommand{\pia}[1]{\pi_{{#1}, T}}
\newcommand{\pimsta}[1]{\pi_{-s, t, T}}

\newcommand{\vXi}[1]{Y_{#1}}

\def\vp{\varphi}
\def\to{\uparrow}
\def\Is{I^{\ms s}}
\def\If{I^{\ms f}}
\def\th{\theta}
\def\pii{\pi}
\def\thi{\th}

\def\su{\subseteq}
\def\De{\Delta}
\def\Xs{X^{\ms s}}
\def\Xss{X^{\ms{bd}}}

\def\Xfa{X^{\ms f, t, T}}
\def\Xsam{X^{\ms s, T}}
\def\Xbam{X^{\ms{bd}, T}}

\def\Xmsta{X^{-s, t, T}}
\def\Xsta{X^{s, t, T}}
\def\Xsaw{X^{\ms {SA}, T}}

\def\Var{\ms{Var}}
\def\Cov{\ms{Cov}}

\def\been{\begin{enumerate}}
\def\enen{\end{enumerate}}
\def\im{\item}
	\def\las{\la^{\ms s}}
	\def\laf{\la^{\ms f}}
	\def\mubd{\nu_{\ms{bd}}^T}
	\def\bec{\begin{corollary}}
	\def\enc{\end{corollary}}
\def\Ts{T^{*}}
\def\Vs{V^{*}}
\def\tff{T \to \ff}

\def\CC{\mathcal{C}}

\def\lrsa{\leftrightsquigarrow}
\def\la{\lambda}

\def\G{\Gamma}
\def\a{\alpha}
\def\s{\sigma}
\def\su{\subseteq}
\def\e{\varepsilon}
\def\t{\tau}
\def\g{\gamma}

\def\de{\delta}

\def\Om{\Omega}
\def\es{\emptyset}
\def\one{\mathbbmss{1}}

\def\De{\Delta}

\def\co{\colon}
\def\ff{\infty}

\def\vp{\varphi}

\def\d{{\rm d}}
\def\k{\kappa}

\def\FF{\mathcal{F}}

\def\f{\frac}

\def\xrd{\xrightarrow{\hspace{.2cm}\ms D\hspace{.2cm}}}

\def\ua{\uparrow}

\def\en{\end}
\def\im{\item}
\def\sm{\setminus}
\def\th{\theta}
\def\bep{\begin{proof}}
\def\enp{\end{proof}}
\def\bepr{\begin{proposition}}
\def\enpr{\end{proposition}}
\def\bec{\begin{corollary}}
\def\enc{\end{corollary}}
\def\bea{\begin{align}}
\newcommand\eea{\end{align}}
\def\beas{\begin{align*}}
\def\eeas{\end{align*}}
\def\bet{\begin{theorem}}
\def\ent{\end{theorem}}
\def\bee{\begin{example}}
\def\ene{\end{example}}

\def\da{\downarrow}
\def\bede{\begin{definition}}
\def\ende{\end{definition}}
\def\bel{\begin{lemma}}
\def\enl{\end{lemma}}
\def\been{\begin{enumerate}}
\def\enen{\end{enumerate}}

\def\beit{\begin{itemize}}
\def\enit{\end{itemize}}
\def\befr{\begin{frame}}
\def\enfr{\end{frame}}

\def\Var{\ms{Var}}
\def\Cov{\ms{Cov}}
\def\ba{\,|\,}

\def\ms{\mathsf}
\def\mb{\mathbb}

\def\Emi{E_k^{M, \ms S}}
\def\Var{\ms{Var}}
\newcommand{\Ei}[1]{E_{{#1, T}}^{\infty, \ms S}}
\newcommand{\Emmi}[1]{E_{#1, T}^{M, \ms S}}
\def\Cov{\ms{Cov}}
\newcommand{\Ec}[1]{E_{{#1, T}}^{\ms{con}}}
\def\CCm{\mathcal C_M}
\newcommand{\Ekm}[1]{E_{{#1, T}}^M}

\def\tktb{\t^{\ms c}_T}
\def\tkktb{\t^{\ms c}_{k, T}}
\newcommand{\vXu}[1]{X_{ #1}}
\newcommand{\vXut}{X_0}
\def\ekmf{E_k^{M \setminus \ff}}

\def\vWui{W_{\ms R}^{T, M}}
\def\Eii{E_k^{\ff, \ms S}}
\def\G{\Gamma}

\def\vli{\lambda_{\ms S}}

\def\d{{\rm d}}
\def\g{\gamma}
\def\P{\mb{P}}
\def\es{\emptyset}
\def\s{\sigma}
\def\la{\lambda}
\def\a{\alpha}

\def\vXui{X^{\ms R}}

\def\eq{\begin{equation}}
\def\en{\end{equation}}

\def\to{\ua}
\def\itf{[0, 1] \rightarrow [0, 1]}

\def\tp{\t^{\ms p}}

\def\If{I^{\ms f}}

\keywords{Continuum percolation,  multi-scale model,  scaling limit,  bounded-hop percolation, wireless communication network}

\begin{document}

\author{Christian Hirsch}
\address[Christian Hirsch]{Bernoulli Institute, University of Groningen, Nijenborgh 9, 9747 AG Groningen, The Netherlands}
\email{christian@math.aau.dk} 
\author{Benedikt Jahnel}
\address[Benedikt Jahnel]{Weierstrass Institute for Applied Analysis and Stochastics, Mohrenstra\ss e 39, 10117 Berlin, Germany}
              \email{benedikt.jahnel@wias-berlin.de} 
\author{Elie Cali}
\address[Elie Cali]{Orange SA, 44 Avenue de la R\'epublique, 92326 Ch\^atillon, France}
              \email{elie.cali@orange.com} 

\title{Percolation and connection times in multi-scale dynamic networks}

\date{\today}

\begin{abstract}
We study the effects of mobility on two crucial characteristics in multi-scale dynamic networks: percolation and connection times. Our analysis provides insights into the question, to what extent long-time averages are well-approximated by the expected values of the corresponding quantities, i.e., the percolation and connection probabilities. In particular, we show that in multi-scale models, strong random effects may persist in the limit. Depending on the precise model choice, these may take the form of a spatial birth-death process or a Brownian motion. Despite the variety of structures that appear in the limit, we show that they can be tackled in a common framework with the potential to be applicable more generally in order to identify limits in dynamic spatial network models going beyond the examples considered in the present work.
\end{abstract}

\maketitle

\section{Introduction}
\label{int_sec}
In  \cite{gilbert}, the \emph{Gilbert graph} was introduced as a network on nodes contained in a metric space, where edges are put between any two nodes closer than a fixed distance. Gilbert thought of this as a toy model for a communication network that helps network operators to study central questions such as: 
\begin{enumerate}
	\item What proportion of nodes is contained in a giant communicating cluster?
	\item What proportion of nodes can connect to a base station in a bounded number of hops?
\end{enumerate}
	The mathematical treatment of these questions falls into the domain of \emph{continuum percolation} \cite{cPerc}.
Despite this early pioneering work, the full potential of spatial random networks in the domain of wireless communication would unfold only  much later.  
The particular traits in communication networks have inspired novel percolation models that are both rooted in applications, and bring also new mathematical facets \cite{sinrPerc,coxPerc}.

											However, one aspect of modern wireless networks has so far received surprisingly little attention: \emph{mobility}. One of the motivations for this line of research is the landmark paper \cite{grossTse}, which showed in an information-theoretic context that mobility can dramatically increase the network capacity. Loosely speaking, if a device starts in a location with bad connectivity, random movements are sufficient to escape this exceptionally bad region after some time.

This breakthrough triggered further research on connection times in large mobile networks with a strong random component such as the analysis in \cite{adHoc} on the average connection times between two devices moving in a fixed domain according to a random waypoint model \cite{bettHart}. Moreover, for devices following a Brownian motion, the first time that the typical device connects to the infinite cluster was investigated in \cite{sousi}. A related Brownian percolation model is also considered in \cite{poisat}. For finite but large random geometric graphs estimates on connection and disconnection times under random walk mobility are provided in \cite{diMi}. However, all of these mathematical works exhibit major shortcomings when viewed from the perspective of real networks, which exhibit \emph{multi-scale features} in a variety of different forms.

\begin{figure}[!htpb]
	\centering
	        \begin{tikzpicture}[scale = 2.1]

	\draw (-1, -1) rectangle (1, 1);

	\fill (-.05,0) -- (0,-.05) -- (.05, 0) -- (0, .05) -- cycle;
	\coordinate[label=-90:{\scriptsize{$o$}}] (A) at (0,0);

	\draw[dashed, ->] (0.2, -.3) -- (1.2, 0.4);
	\fill[green] (0.2, -.3) circle (1pt);
	\draw (0.2, -.3) circle (1pt);

	\draw[dashed, ->] (0.2, .8) -- (-.6, -1.1);
	\fill[green] (0.2, .8) circle (1pt);
	\draw (0.2, .8) circle (1pt);

	\draw[dashed, ->] (-.7, -.4) -- (.6, -1.1);
	\fill[green] (-0.7, -.4) circle (1pt);
	\draw (-0.7, -.4) circle (1pt);


	\draw[dashed, ->] (-.55, .35)--(.5, .2);
	\draw[dashed, ->] (.85, .45)--(.5, .3);

	\fill[cyan] (0.8, .4) rectangle (0.9, .5);
	\draw (0.8, .4) rectangle (0.9, .5);
	\fill[cyan] (-.6, .3) rectangle (-.5, .4);
	\draw (-.6, .3) rectangle (-.5, .4);
	
\end{tikzpicture}
\hspace{.1cm}
\begin{tikzpicture}[scale = 2.1]

	\draw (-1, -1) rectangle (1, 1);

	\fill (-.05,0) -- (0,-.05) -- (.05, 0) -- (0, .05) -- cycle;
	\coordinate[label=-90:{\scriptsize{$o$}}] (A) at (0,0);


	\draw[dashed, ->] (0.8, -.1) -- (-.6, -1.1);
	\fill[green] (0.8, -.1) circle (1pt);
	\draw (0.8, -.1) circle (1pt);

	\draw[dashed, ->] (.7, -.4) -- (-1.1, .5);
	\fill[green] (0.7, -.4) circle (1pt);
	\draw (0.7, -.4) circle (1pt);


	\draw[dashed, ->] (-.25, -.35) -- (-.9,-.7);
	\fill[cyan] (-.3, -.4) rectangle (-.2, -.3);
	\draw (-.3, -.4) rectangle (-.2, -.3);

	\draw[dashed, ->] (-.55, .35)--(.5, .2);
	\fill[cyan] (-.1, .23) rectangle (0, .33);
	\draw (-.1, .23) rectangle (0, .33);

	\draw[dashed, ->] (.85, .45)--(.5, .3);
	\fill[cyan] (0.55, .28) rectangle (0.65, .38);
	\draw (0.55, .28) rectangle (0.65, .38);

\end{tikzpicture}
\hspace{.1cm}
\begin{tikzpicture}[scale = 2.1]

	\draw[white,dashed, ->] (.7, -.4) -- (-1.1, -1.1);
	\draw (-1, -1) rectangle (1, 1);

	\fill (-.05,0) -- (0,-.05) -- (.05, 0) -- (0, .05) -- cycle;
	\coordinate[label=-90:{\scriptsize{$o$}}] (A) at (0,0);

	\draw[dashed, ->] (.7, -.8) -- (1.1, .5);
	\fill[green] (0.7, -.8) circle (1pt);
	\draw (0.7, -.8) circle (1pt);

\draw[dashed, ->] (-.25, -.35) -- (-.9,-.7);
	\fill[cyan] (-.6, -.52) rectangle (-.7, -.62);
	\draw (-.6, -.52) rectangle (-.7, -.62);

	\draw[dashed, ->] (-.55, .35)--(.5, .2);
	\fill[cyan] (.35, .16) rectangle (.45, .26);
	\draw (.35, .16) rectangle (.45, .26);

\end{tikzpicture}
		       \caption{Two-scale mobility model. Slow users (blue squares) together with passing fast users (green circles). }

		        \label{2scale_fig}
\end{figure}
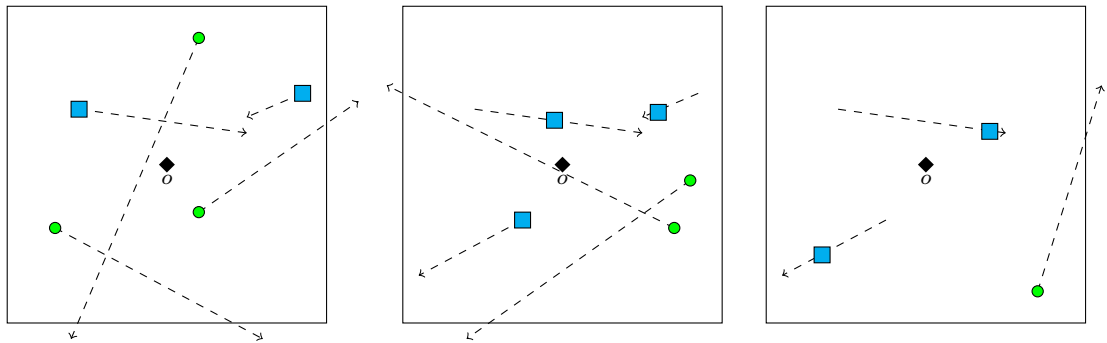

As a first example, note that the literature described in the preceding paragraph all assume a single movement model over the entire time horizon. However, this assumption clearly does not mirror a typical working day, beginning with a commute over a long distance followed by predominantly local movements, and finally concluded again by a long commute. We address this shortcoming by introducing a \emph{two-scale mobility model} capturing the type of such more realistic movement patterns, as illustrated in Figure \ref{2scale_fig}. 

\begin{figure}[!htpb]
	\centering
	        \begin{tikzpicture}[scale = 1.4]
\fill (0, 0) circle (1pt);
\draw (0, 0)--(0.03, -0.05)--(0.00, -0.04)--(-0.04, -0.02)--(-0.08, 0.02)--(-0.03, 0.02)--(-0.02, 0.05)--(0.01, 0.04)--(0.03, 0.04)--(0.08, 0.05)--(0.06, 0.09)--(0.05, 0.05)--(0.03, 0.04)--(0.06, 0.06)--(0.03, 0.09)--(0.06, 0.08)--(0.06, 0.11)--(0.05, 0.13)--(0.01, 0.13)--(0.00, 0.10)--(0.04, 0.06)--(0.03, 0.01)--(0.03, 0.05)--(0.02, 0.09)--(-0.02, 0.05)--(0.01, 0.07)--(0.03, 0.07)--(0.03, 0.11)--(0.02, 0.12)--(0.04, 0.11)--(0.05, 0.12)--(0.04, 0.16)--(0.07, 0.20)--(0.11, 0.18)--(0.13, 0.14)--(0.10, 0.16)--(0.12, 0.14)--(0.17, 0.09)--(0.15, 0.06)--(0.15, 0.10)--(0.13, 0.10)--(0.12, 0.13)--(0.17, 0.13)--(0.17, 0.16)--(0.21, 0.12)--(0.19, 0.14)--(0.16, 0.14)--(0.21, 0.19)--(0.23, 0.19)--(0.23, 0.15)--(0.28, 0.19);
\draw[dashed] (0.03, -0.05) circle (2cm);
\draw[dashed] (0.00, -0.04) circle (2cm);
\draw[dashed] (-0.04, -0.02) circle (2cm);
\draw[dashed] (-0.08, 0.02) circle (2cm);
\draw[dashed] (-0.03, 0.02) circle (2cm);
\draw[dashed] (-0.02, 0.05) circle (2cm);
\draw[dashed] (0.01, 0.04) circle (2cm);
\draw[dashed] (0.03, 0.04) circle (2cm);
\draw[dashed] (0.08, 0.05) circle (2cm);
\draw[dashed] (0.06, 0.09) circle (2cm);
\draw[dashed] (0.05, 0.05) circle (2cm);
\draw[dashed] (0.03, 0.04) circle (2cm);
\draw[dashed] (0.06, 0.06) circle (2cm);
\draw[dashed] (0.03, 0.09) circle (2cm);
\draw[dashed] (0.06, 0.08) circle (2cm);
\draw[dashed] (0.06, 0.11) circle (2cm);
\draw[dashed] (0.05, 0.13) circle (2cm);
\draw[dashed] (0.01, 0.13) circle (2cm);
\draw[dashed] (0.00, 0.10) circle (2cm);
\draw[dashed] (0.04, 0.06) circle (2cm);
\draw[dashed] (0.03, 0.01) circle (2cm);
\draw[dashed] (0.03, 0.05) circle (2cm);
\draw[dashed] (0.02, 0.09) circle (2cm);
\draw[dashed] (-0.02, 0.05) circle (2cm);
\draw[dashed] (0.01, 0.07) circle (2cm);
\draw[dashed] (0.03, 0.07) circle (2cm);
\draw[dashed] (0.03, 0.11) circle (2cm);
\draw[dashed] (0.02, 0.12) circle (2cm);
\draw[dashed] (0.04, 0.11) circle (2cm);
\draw[dashed] (0.05, 0.12) circle (2cm);
\draw[dashed] (0.04, 0.16) circle (2cm);
\draw[dashed] (0.07, 0.20) circle (2cm);
\draw[dashed] (0.11, 0.18) circle (2cm);
\draw[dashed] (0.13, 0.14) circle (2cm);
\draw[dashed] (0.10, 0.16) circle (2cm);
\draw[dashed] (0.12, 0.14) circle (2cm);
\draw[dashed] (0.17, 0.09) circle (2cm);
\draw[dashed] (0.15, 0.06) circle (2cm);
\draw[dashed] (0.15, 0.10) circle (2cm);
\draw[dashed] (0.13, 0.10) circle (2cm);
\draw[dashed] (0.12, 0.13) circle (2cm);
\draw[dashed] (0.17, 0.13) circle (2cm);
\draw[dashed] (0.17, 0.16) circle (2cm);
\draw[dashed] (0.21, 0.12) circle (2cm);
\draw[dashed] (0.19, 0.14) circle (2cm);
\draw[dashed] (0.16, 0.14) circle (2cm);
\draw[dashed] (0.21, 0.19) circle (2cm);
\draw[dashed] (0.23, 0.19) circle (2cm);
\draw[dashed] (0.23, 0.15) circle (2cm);
\draw[dashed] (0.28, 0.19) circle (2cm);
\fill[red] (2.75, 2.43) rectangle (2.85, 2.53);
\fill[red] (0.01, 1.04) rectangle (0.11, 1.14);
\fill[red] (0.79, -0.11) rectangle (0.89, -0.01);
\fill[red] (-3.38, 1.94) rectangle (-3.28, 2.04);
\end{tikzpicture}\qquad
	        \input{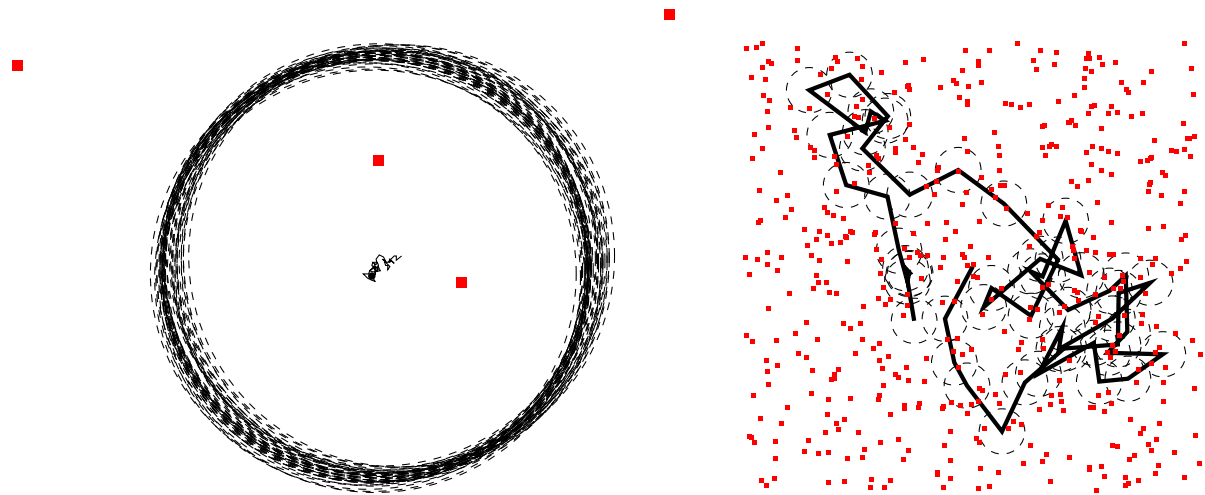}
		       \caption{Infrastructure-augmented model. Typical moving node (black line) inside process of sparse (left) or dense (right) sinks (red squares). Dashed circles indicate communication range via relaying. For clarity, relays themselves are not shown.
		       }
		        \label{inf_fig}
\end{figure}

As a second example, it is not realistic to assume that a large-scale wireless network can rely entirely on the form in which data forwarding is considered in continuum percolation. In practice, we will witness the relaying technology being employed inside an \emph{infrastructure-augmented network}. Thus, starting from a randomly scattered configuration of base stations or sinks, nodes may rely on other nodes as relays to establish communication with a base station. However, massive use of relaying also comes at a price: Each hop induces a certain jitter and delay, so that in practice, network operators will put a constraint on the admissible number of hops. When considering large sink-augmented networks over a long time, then both the number of base stations as well as the time horizon grow, but possibly at completely different scales, as illustrated in Figure \ref{inf_fig}. For clarity, we here approximate the communication ranges via relaying by disks of some radius. Hence, we witness a further incarnation of a multi-scale phenomenon.

We will lay the groundwork for a rigorous asymptotic analysis over long time horizons of two key performance indicators of large dynamic communication networks: percolation and connection times. Loosely speaking, the \emph{percolation times} are the points in time when a typical node connects to a positive proportion of the other network participants via relay communication; the \emph{connection times} are the points in time where a typical node can connect to some static sink via relay communication.

Often, it is possible to compute, or at least approximate, how a wireless network behaves in expectation. A commonly held belief is that, when averaged over a long time horizon,  the key network characteristics in dynamic random network are represented fairly accurately through their expected values. We challenge this view and draw a more nuanced picture through an asymptotic analysis. While to a certain extent, the time averaging removes some random fluctuations, we identify clearly specific sources of randomness that are so pronounced that they persist over time. 

The rest of the manuscript is organized as follows.  In Section \ref{mod_sec}, we introduce precisely the two-scale mobility model and the infrastructure-augmented network. We also state Theorem \ref{thm_2} and \ref{thm_1} as the main results in the present work. In Section \ref{sec_outline}, we outline the proof based on a common general blueprint to tackle weak convergence problems in spatial random networks. Finally, Section \ref{sec_proofs} contains the detailed proofs.

\section{Setting and main results}
\label{mod_sec}
First, in Section \ref{sys_sec}, we introduce the general stochastic-geometry framework which forms the basis of our investigations. Then, in Sections \ref{perc_sec} and \ref{khop_sec}, we describe the specific two-scale mobility and infrastructure-augmented models for which we perform the long-time analysis.

%
%
\subsection{General system model} 
\label{sys_sec}
The network nodes are scattered at random across the entire Euclidean plane. Mathematically speaking, they form a homogeneous Poisson point process $X = \{X_i\}_{i\ge 1}$ with some intensity $\la > 0$. In the present system model, nodes move independently, which is captured by attaching to each $X_i$ a trajectory $\G_i\co [0, T] \rightarrow \R^d$ taken from an i.i.d.~collection of paths $\{\G_i\}_{i \ge 1}$
$$X_i(t) := X_i + \G_i(t).$$ 
More precisely, the trajectories are i.i.d.~copies of the random element $\G = \{\G(t)\}_{t \le T}$ in the space $\M$ of c\`adl\`ag (right-continuous with left limits) paths from $[0, T]$ to $\R^d$ starting at the origin. The space $\M$ is a Polish space when endowed with the Skorokhod topology and serves as the mark space for the independently marked homogeneous Poisson point process. At this point, we keep the mobility model still general, and specify it further below.

 We think of a \emph{typical node} as a distinguished node starting the movement at the origin, i.e., $X_0 = o$, and note that this amounts to considering the process under its Palm distribution. The relay connection is modeled through a simple \emph{distance-based criterion}: $X_i$ and $X_j$ can communicate directly if their Euclidean distance is smaller than some connection threshold $r > 0$, i.e., two nodes are connected at time $t \le T$ if $|X_i(t) - X_j(t)| \le r$. Additionally, we allow communication through relaying in several hops. This form of connectivity is encoded in the Gilbert graph, where an edge is drawn between any two nodes at distance at most $r$.

%
%
\subsection{Percolation times in the two-scale mobility model}
\label{perc_sec}

In the vast majority of the literature, a single mobility model is assumed throughout the entire time horizon. However, this is not what we experience in everyday life, where mobility occurs at several scales. For instance, on a normal work day, nodes may move over long distances when commuting to work, then are close to being static within their offices, and finally move again over long distances on the commute home. We think of this as an alternation between two phases: a \emph{fast} one, corresponding to a long distance to run through, and a \emph{slow} one, corresponding to a short distance. We incorporate this observation into the mobility model by relying on a sequence of waypoints that each node visits successively by traveling along the connecting line segments. The two-scale nature of the mobility is implemented by imposing that the node alternates between moving to distant waypoints and to local ones that  are in the vicinity of the current location. 

More precisely, for each node $X_i$, let $\{V_{i,j}\}_{j \ge 0}$ denote an i.i.d.~sequence of displacement vectors drawn from an isotropic probability distribution $\k(\d v)$ of bounded support, which is absolutely continuous with respect to the Lebesgue measure on $\R^d$. Then, we use the total time horizon $T$ as scaling parameter for the two-scale mobility model. That is, we obtain the induced sequences of waypoints $ (P_{i,0}, P_{i,1}, P_{i,2}, P_{i,3}, \dots),$ where $P_{i,0}  :=  X_i$, and 
$$P_{i,j + 1} :=  P_{i,j} + T^{ j\, \ms{mod}\, 2} V_{i,j}.$$
Without loss of generality, we will fix the high speed at 1, and take for the slow speed the same ratio as for the distances between the waypoints.
{Now, if $T_{i, j}^0 := T_{i, j}/T$ denotes the rescaled arrival time at $P_{i, j}$, then $T_{i, 0} = 0$ and }
$$ \frac{|P_{i, j+1} - P_{i, j}|}{T_{i, j+1}^0 - T_{i, j}^0} = T^{ j\, \ms{mod}\, 2}, $$
 i.e., $ (T_{i,0}, T_{i,1}, T_{i,2}, T_{i,3}, \dots) = (TT_{i,0}^0, TT_{i,1}^0, TT_{i,2}^0, TT_{i,3}^0, \dots)$, where  
\begin{align*}
	T_{i,j + 1}^0 :=  T_{i,j}^0 +  |V_{i,j}|.
\end{align*}
{To summarize, we link the length of the time horizon $T$ to the waypoint-scaling in a way that on average there are a finite and non-vanishing number of waypoints within the time horizon.}

Note that this model could be generalized, for instance by assuming that the waypoints are constructed as $P_{i, j+1} = P_{i, j} + K_1(T)^{j \,\ms{mod}\, 2} V_{i, j}$, whereas the speed in phase $j$ equals $v_j = K_2(T)^{j\, \ms{mod}\, 2} v$ for some speed $v > 0$ and scalings $K_1(T), K_2(T)$. For the sake of simplicity, we restrict our attention to the most elementary case where   $K_1(T) = K_2(T) = T$, and without loss of generality we suppose $v = 1$.
Then, the path $\G_i$ of node $X_i$ consists of straight lines between the sequence of its waypoints, alternating between nearby and very distant targets, and thereby constituting the two-scale nature of the mobility model. In particular, the arrival times $\{T_{i,j}\}_{j \ge 0}$ form a renewal process with inter-arrival distribution $T|V|$ where $V \sim \k(\d v)$. The randomness enters the system only via the initial positions $\{X_i\}_{i \ge 1}$ and the displacement vectors $\{V_{i, j}\}_{j \ge 0}$.

 Now, we define the \emph{percolation time}
 \begin{align}      
 \big|\{t \le T:\,  o\lrsa_t \ff\}\big| \end{align} 
	 \noindent in a horizon $[0, T]$ as the total amount of time that the origin is contained in an unbounded connected component. More precisely, $o \lrsa_t \ff$ means that $o$ is contained in an unbounded component of the Gilbert graph formed on $\{o\} \cup \{X_i(t)\}_{i \ge 1}$. To ease notation, we let $\pii(X(t))$ denote the indicator of the event that at time $t$, starting from the static origin there exists an unbounded sequence of connected nodes in $X(t)$.
%
%
In fact our methods allow to treat a measure-valued refinement of the connection time, namely the \emph{empirical percolation time measure}
$$\tp_T :=\frac 1T \int_0^T\pii(X(t))\de_{t/T} \d t= \int_0^1\pii(X(tT))\de_t \d t, $$
where $\de_t$ denotes the Dirac measure in $t\in [0,1]$. In words, $\tp_T$ is the measure of (rescaled) times that the origin is contained in the unbounded connected component of the node process. We consider $\tp_T$ as a random element in the space of measures on $[0, 1]$, which will be equipped with the topology of weak convergence. This allows us, for example, to evaluate the percolation time of the typical node via $T\tp_T([0,1])$.

Now, we revisit the motivating question:
\begin{center}
	        \emph{Is the time-averaged percolation time approximated well by the percolation probability?}
\end{center}
In other words, do we see a law-of-large-numbers type of situation in which the time average approaches the expected value? During the fast phase, this is highly plausible since there are substantial portions of times where the typical node sees a completely new environment, thereby leading to rapid decorrelation. However, during the slow phase, the typical node moves inside a relatively restricted spatial domain, so that the family of potentially relevant relays will exhibit substantial overlaps at different points in time. This creates high correlations persisting in the average over time. The limiting time-averaged percolation time is described through a spatial birth-and-death process corresponding to nodes entering and exiting their slow phase.

In the main result for this model, Theorem \ref{thm_2} below, we determine the empirical percolation time measure $\tp_T$ in the long-time limit, i.e., $T \to \ff$. In this limit, it seems that once nodes switch into a slow phase, they appear at a certain location out-of-nowhere, move at constant speed one, and after the phase is completed, they disappear immediately. In addition, the nodes in the global phase form an ephemeral background process of points that may be observed for a single instant, but are afterwards never seen again.

%
%
To make this picture precise, we first describe the limit during the slow phases. Let $\Xss=\{(\Xss_i,\s_i,V_i)\}_{i \ge 1}$ be a Poisson point process in $\R^d\times [0, 1]\times \R^d$ with intensity measure\linebreak $\nu([0, 1])^{-1}\la\d x\otimes\nu(\d t)\otimes\k(\d v)$, where 
$$\nu(\d t) := \E\Big[\sum_{j \ge 0}\de_{T_{0,2j}^0}(\d t)\Big]$$ 
denotes the intensity measure of a renewal process on $[0,\infty)$ with inter-arrival times distributed according to $V'+ V''$, with $V'$ and $V''$ being i.i.d.~with distribution given by $\k(\d v)$. 
In words, this Poisson point process represents a birth-and-death process with locations $\Xss_i$ and arrival times $\s_i$ such that after time $\s_i$ a particle is born at $\Xss_i$ and then moves to location $\Xss_i+V_i$ with speed one. The convolution comes from the alternation between slow and fast phases.  
Let 
$$\Xs(t): = \{\Xss_i + (t - \s_i)V_i/|V_i|\}_{i\co \s_i\le t\le \s_i +|V_i|}$$ 
denote the process of nodes in the limiting process of slow nodes at time $t$. By the displacement theorem \cite[Exercise 5.1]{poisBook}, $\Xs(t)$ is a homogeneous Poisson point process with $t$-dependent intensity 
$\las(t) := \la \P(t \in \Is), $
where 
$\Is := \cup_{j \ge 0}[T_{0,2j}^0, T_{0,2j + 1}^0]$
denotes the slow phase of a typical node.

%
%
During the fast phases the nodes move so far that they achieve perfect decorrelation. Hence, they can be treated as a homogeneous Poisson point process that is marginalized at every time instant. More precisely, 
$$\thi(\Xs(t); \laf(t)) := \E[\pii(\Xs(t) \cup X') \ba \Xs(t)]$$
denotes the conditional percolation probability at the origin in the process $\Xs(t) \cup X'$ where $X'$ is an independent homogeneous Poisson point process with intensity 
$\laf(t) := \la \P(t \in \If), $
where 
$\If := \cup_{j \ge 0}[T_{0,2j + 1}^0, T_{0,2j + 2}^0]$
denotes the fast phase of a typical node.
In particular,  $\thi(\la) := \E[\thi(\Xs(t); \laf(t))]$
does not depend on $t$ and recovers the classical percolation probability from continuum percolation \cite{cPerc}. Let $\la_c:=\inf\{\la'\ge 0\colon\thi(\la') > 0\}$ denote the critical intensity for percolation. We now describe the asymptotic empirical percolation time measure $\tp_T$ as $T \ua \ff$, where we write $\xrd$ for convergence in distribution in the limit $T \to \ff$.
	
%
%
\begin{theorem}[Asymptotic percolation time]\label{thm_2}
	Let  $\la > \la_c$. Then, as $T \to \ff$,
$$
	\tp_T\xrd  \thi(\Xs(t);{\laf(t)}) \d t.
$$
\end{theorem}
We present the proof of Theorem~\ref{thm_2} in Section \ref{sec_outline} based on a sequence of key proposition, that are proved in Section \ref{sec_proofs}.

%
%
%
\subsection{Connection times in the infrastructure-augmented model}
\label{khop_sec}
Despite the rise of relaying technology, also future communication networks will still rely on some form of infrastructure as the sink locations for data transmission. Hence, the percolation time studied in Section \ref{perc_sec} is essentially a proxy for the more realistic connection time to some sink. In this manuscript, we assume that the sinks also form a homogeneous Poisson point process $Y = \{Y_j\}_{j \ge 1}$ with some intensity $\vli > 0$.

To compensate for the additional complexities of the connection time, we restrict the mobility model to a standard continuous-time random walk. That is, nodes choose random waypoints sequentially according to the probability measure $\k(\d v)$ and directly jump to them after exponentially distributed waiting times. We assume that the trace of the coordinate-covariance matrix associated with random vectors from $\k(\d v)$ equals $d$. This normalization will later ensure convergence to a standard Brownian motion.

As indicated in Section \ref{int_sec}, the investigation of large dynamic sink-augmented networks over a long time yields a second class of prototypical examples for multi-scale networks, since both the time horizon as well as the density of sinks grow, but possibly at completely different scales.

%
%

 We assume that sinks can communicate with nodes directly within the connection threshold $r$.
 Additionally, nodes can connect to sinks via relaying through other nodes. However, in order to avoid long delays and jitter in the data transmission, network operators typically put constraints on the number of admissible hops. More precisely, $x \stackrel k{\lrsa_t} y$ means that $x$ is at graph distance at most $k$ from $y$ in the Gilbert graph formed on $\{x,y\} \cup \{X_i(t)\}_{i \ge 1}$. Then, we define the \emph{connection time}
\begin{align} \big|\{t \le T:\, X_0(t)\stackrel k{\lrsa_t} Y_j \text{ for some $Y_j \in Y$}\}\big|, \end{align}
	\noindent in a time horizon $[0, T]$, as the total amount of time that the typical node can connect in at most $k \ge 1$ relaying hops to at least one sink $Y_j \in Y$. 
%
%
%
%
In fact, our methods allow again to describe the more refined \emph{empirical $k$-hop connection measure} 
$$
\tkktb : = \frac1T \int_0^T\one\{\vXut(t)\in\Xi^k(t)\}\de_{t/T} \d t = \int_0^1\one\{\vXut(tT)\in\Xi^k(tT)\}\de_t \d t$$
over the time horizon $[0, T]$, where 
$$\Xi^k(t) : = \{x \in \R^d:\, x \stackrel k{\lrsa_t} Y_j\text{ for some }Y_j\in Y\}$$
denotes the set of all points connecting to a sink in at most $k$ hops at time $t$. 

Getting insights into the $k$-hop connection time is essential for network operators. Indeed, while it is unrealistic to exclude times of bad connectivity entirely, at least on average over a long time, the system should offer excellent quality of service to its nodes. Since the system is time-stationary, the expected value of $\tkktb$ recovers the \emph{$k$-hop connection probability}
\begin{align}
	\label{th_eq}
	 \E[\tkktb([0,1])] = \P(o \in\Xi^k(0)).
\end{align}
%
%

%
%
However, since the empirical $k$-hop connection measure $\tkktb$ depends sensitively on all model parameters in a complex manner, a closed-form expression of its distribution is out of reach. In the present work, we elucidate how to obtain a conceptually clean and readily interpretable description in a scenario where the sink density is relatively low, so that the system critically relies on the option of relaying to achieve good quality of service. More precisely, we assume that $k \ua \ff$ and $\vli \da 0$ such that
\begin{align}
	\label{scale_la_eq}
	\vli k^d = c_0,
\end{align}
for some fixed $c_0 > 0$. Now, we again revisit the motivating question:
\begin{center}
	    \emph{Is the time-averaged connection time approximated well by the connection probability?}
\end{center}

We show that the answer to this question depends sensitively on the sparseness of the sinks in relation to the time horizon. If the sinks are sufficiently dense, then the typical node sees a large number of different sinks during the movement, thereby leading to near-perfect decorrelation. On the contrary, for sparsely distributed sinks, the typical node basically sees the same configuration during the entire time horizon. Finally, in a particularly subtle critical scaling, the typical node encounters a finite, random number of novel sinks, which again leads to substantial correlations that are captured in the limit by a Brownian motion navigating through a Poisson process of sinks.

%
%
Henceforth, $B_r(x) : = \{y \in \R^d:\, |y - x| \le r\}$ denotes the Euclidean ball with radius $r > 0$ centered at $x \in \R^d$, and we write $B_r:=B_r(o)$. We now elucidate in greater detail, how the scaling \eqref{scale_la_eq} allows to reduce the intensity of sinks by relying more aggressively on forwarding data. To that end, we assume that the process $\vXu{}$ of nodes is in the super-critical regime of percolation \cite{cPerc}, i.e., $\la>\la_c$. That is, at any fixed time $t \le T$, with probability 1, the Boolean model $\bigcup_{i\ge 1}B_{r/2}(\vXu i(t))$ contains a unique unbounded connected component of nodes $\CC(t) \su \vXu{}(t)$.

%
A key feature of the super-critical percolation phase is the shape theorem \cite[Theorem 2]{yao}. Loosely speaking, it states that the number of hops needed to connect two nodes in $\CC(t)$ grows linearly in the Euclidean distance. More precisely, for nodes $\vXu i(t)$, $\vXu j(t)$ that are in the infinite connected component $\CC(t)$, we let $T(\vXu i(t), \vXu j(t))$ denote the minimum number of hops of distance at most $r$ needed to move from $\vXu i(t)$ to $\vXu j(t)$. Then, $T(\cdot, \cdot)$ scales asymptotically linearly in the distance of the two nodes. More precisely, there exists a deterministic \emph{stretch factor} $\mu > 0$ such that almost surely
\begin{align}
	\label{stretch_eq}
	\lim_{\substack{|\vXu i(t) - \vXu j(t)| \ua \ff \\ \vXu i(t), \vXu j(t) \in \CC(t)}} \frac{T(\vXu i(t) , \vXu j(t))}{|\vXu i(t) - \vXu j(t)|} =  \mu.
\end{align}
%
%
The shape theorem unveils why \eqref{scale_la_eq} is the correct scaling to balance the sparsity of sinks and the number of hops. Indeed, it implies that with high probability, asymptotically, the locations that can be reached with at most $k$ hops from a point in the unbounded connected component are contained in a $k/\mu$-ball around this point. Then, \eqref{scale_la_eq} encodes that on average this ball contains a constant number of sinks.

%
%
In the super-critical regime, \cite[Theorem 2.6]{bhopPerc} captures the $k$-hop connection probability. In a nutshell, a $k$-hop connection is possible if both the typical node and at least one sink in the $k/\mu$-ball are in the unbounded connected component. Leveraging the void probability of the Poisson point process shows that if $k \ua \ff$ and $\vli \da 0$ such that the scaling \eqref{scale_la_eq} holds, then
$$\lim_{k \to \ff} \P(o \in\Xi^k(0)) = \th(\la) \big(1 - \exp(- c_0\th(\la) |B_{1/\mu}|)\big).$$

%
%
Although this scaling describes in closed form the asymptotic $k$-hop connection probability, relying on it in the setting of mobile nodes may be misleading. Indeed, although \eqref{th_eq} shows that the $k$-hop connection probability is the expected average $k$-hop connection time, this picture can be severely inaccurate over short time horizons. Indeed, then the averaging nature of the mobility may not be visible effectively, so that random fluctuations dominate the picture. 
To work out the impact of mobility cleanly, we need to take into account the sink-densities relation to the considered time horizon. More precisely, we investigate scalings of the form
\begin{align}
	\label{scale_t_eq}
	\vli(T) = T^{-\a},
\end{align}
for some parameter $\a > 0$, controlling the sink density. As we will see below in Theorem \ref{thm_1}, the invariance principle of the node movement gives rise to a phase transition with the Gaussian scaling at exponent $\a = d/2$.

%
%
The main result for this model describes the scaling of the asymptotic $k$-hop connection time for different scalings of the sink densities. Here, Theorem \ref{thm_1} unveils accurately how the degree of averaging influences the amount of randomness remaining in the limit. More precisely, when sinks are dense, then all randomness disappears in the sense that the averaged $k$-hop connection time converges to the deterministic connection probability. In contrast, when sinks are sparse, the typical node 
moves so little that the set of relevant sinks stays invariant during the entire horizon, thereby giving rise to a Poisson random variable in the limit. The critical case is right in the middle. Here, the set of relevant sinks does change within the horizon, but only so slowly that some degree of randomness is still visible in the limit.

To ease notation, we write $X'(A) = \#(X' \cap A)$ for the number of points of a point process $X'$ in a Borel set $A \su \R^d$. Under the scalings \eqref{scale_la_eq} and \eqref{scale_t_eq} we also write $\tktb$ instead of $\tkktb$.
\bet[Asymptotic $k$-hop connection times]
\label{thm_1}
Let $\la > \la_c$ and assume the multi-scale regime encoded by \eqref{scale_la_eq} and \eqref{scale_t_eq}.\\[2ex]
{\bf Dense sinks.}
If $\a < d/2$, then, as $T \to \ff$,
  \begin{align}
	  \label{long_eq}
	   \tktb \xrd \th(\la)\big(1 - \exp( - c_0\th(\la)|B_{1/\mu}|)\big) \d t.
  \end{align}
{\bf Sparse sinks.}
If $\a > d/2$, then, as $T \to \ff$,
\begin{align}
	\label{short_eq}
	\tktb \xrd \th(\la)\big(1 - (1 - \th(\la))^N\big) \d t,
\end{align}
where $N$ is a Poisson random variable with parameter $c_0|B_{1/\mu}|$.\\
{\bf Critical density.}
If $\a = d/2$, then, as $T \to \ff$,
  \begin{align}
	  \label{crit_eq}
	  \tktb \xrd \th(\la) \big( 1- (1 - \th(\la))^{Y'(B_{1 / \mu}(W_t))}\big) \d t,
  \end{align}
 where $Y'$ is a homogeneous Poisson point process with intensity $c_0$ and $W_t$ is a standard Brownian motion. 
\ent

In the following section we present the proof of Theorem~\ref{thm_1} based on a sequence of key propositions that are proved in Section \ref{sec_proofs}.

\noindent {\bf Remark.} Many of the arguments in Theorems \ref{thm_2} and \ref{thm_1} generalize also to settings, where the nodes $\vXu{}$ do not necessarily form a Poisson point process. A particularly relevant case from the point of view of wireless networks is that of a \emph{Cox point process}. Here, the nodes represent devices that are constraint to lie on the street system of a large city, which is modeled as a spatial random network \cite{coxPerc}.

\section{Outline}\label{sec_outline}
Although on the surface, Theorems \ref{thm_2} and \ref{thm_1} involve limiting regimes of rather distinct flavors, we show that they can both be tackled under a common umbrella within a three-step strategy. In fact, this strategy holds the promise of serving as a general blueprint for approaching distributional limits appearing in spatial models for wireless networks. The three steps can be succinctly summarized as follows. Here, $\t $ denotes either the time-averaged percolation or the time-averaged connection time.
\been
\im {\bf Identification of persistent information.} In a first step, some intuition into the specific problem at hand is needed to decide what kind of information might persist over long time horizons. This information is gathered in a $\s$-algebra $\FF$.
\im {\bf Conditional decorrelation.} In a second step, one shows that after conditioning on $\FF$, percolation or connection events that are macroscopic time instances apart decorrelate in the limit. Then, a second-moment method reduces the problem to understanding the distributional limit of the conditional expectation $\E[\t \ba \FF ]$.
\im {\bf Computation of conditional expectation.} Finally, the conditional expectation $\E[\t \ba \FF]$ is computed and a concise expression for its limiting distribution is derived.
\enen

In the following, we illustrate how to implement this blueprint for the two-scale mobility and infrastructure-augmented model.

%
%
\subsection{Two-scale mobility model}
\label{sec_pre_sf}
Before starting to implement the above blueprint for the two-scale mobility network, we insert a preliminary step, reducing the proof of Theorem~\ref{thm_2} from percolation to infinity to percolation outside an $M$-box. More precisely, let $\pi^M(X(tT))$ denote the indicator of the event that at time $t$, there exists a sequence of connected nodes in $X(tT)$ that create a path between $o$ and a node at distance at most 1 from $\partial Q_M$, the boundary of the box $Q_M$ with side-length $M>0$, centered at the origin. The quantities $(\tp_T)^M$, $\pi^M$ and $\th^M$ are understood accordingly.
\bepr[Asymptotic percolation time -- $M$-approximation]
\label{Thm_1m} 
Let $M > 0$. Then, as $\tff$,
$$
(\tp_T)^M \xrd \th^M(\Xs(t); \laf(t)) \d t.
$$ 

\enpr
%
%
The proof of the main Proposition~\ref{Thm_1m} will be discussed further below. Let us now use Proposition~\ref{Thm_1m} to prove Theorem~\ref{thm_2}. Henceforth, we simplify notation and let
$$\pia t^M := \pi^M(X(tT))\quad\text{ and }\quad\pia t := \pi(X(tT))$$
denote the percolation and the $M$-percolation indicators at time $tT$.

%
%
  \bep[Proof of Theorem~\ref{thm_2}]
      Let $f$ be a bounded Lipschitz function with Lipschitz constant 1 and  $g\co\itf$ be continuous.
      We first consider the difference on the left-hand side
      \begin{align*}
				\E[|f((\tp_T)^M(g)) - f(\tp_T(g))|]
				\le \int_0^1(\E[\pia t^M] - \E[\pia t])g(t)\d t
				\le \th^M(\la) - \thi(\la),
      \end{align*}
      where we could exchange integrals due to the boundedness of the integrand. Note that the second inequality follows because $X(tT)$ is a homogeneous Poisson point process with intensity $\la$ at every time instant $tT$. Since the expression on the right does not depend on $T$, we obtain convergence in $M$ uniformly over all $T > 0$. 

			Similarly, write $Z^M$ and $Z^\ff$ for the right-hand sides in Proposition \ref{Thm_1m} and Theorem \ref{thm_2}. Then, 
      \begin{align*}
        |\E[f(Z^M) - f(Z^\ff)]| 
				\le\int_0^1 \big(\E[\th^M(\Xs(t);{\laf(t)})] - \E[\th(\Xs(t);{\laf(t)})]\big)g(t)\d t
        \le \th^M(\la) - \thi(\la),
      \end{align*}
      which again tends to 0 as $M \to \ff$.
      \enp

The proof of Proposition~\ref{Thm_1m} follows in three main steps, represented by Propositions \ref{decorr_prop}, \ref{fast_prop} and \ref{slow_prop} below. To ease notation, we henceforth write $\pi$ and $\th$ instead of $\pi^M$ and $\th^M$.

\subsubsection{Identification of persistent information}
Loosely speaking, the intuition is that the persisting information $\FF$ is captured by the nodes entering the $M$-box during their slow phase. A critical ingredient in the following steps is that nodes visit a spatial neighborhood only once with a high probability. 
More precisely, a path $\G$ is \emph{$T$-self-avoiding} if 
$$Q_{4M}(\G(T_{2j})) \cap [\G(T_{2j' - 1}), \G(T_{2j'})] = \es$$
for all $j, j' \ge 0$, with $j' \ne \{j, j + 1\}$ and $T_{2j}, T_{2j' - 1} \le T$, where we suppressed the first index for brevity. That is, during the fast movement, the trajectory never crosses the $4M$-neighborhood of a waypoint in which the slow phase begins. 
Finally, we let 
$$\Xsaw := \{X_i\co \text{$\G_i$ is $T$-self-avoiding}\}$$
denote the family of all \emph{$T$-self-avoiding nodes}. 

The first step is to note that during their fast movement, the nodes have outstanding mixing properties, so that in the limit only the randomness from the slow movements persists. To that end, let
 $$\Ts_i := \inf\{T_{i, 2j}\co X_i(T_{i, 2j}) \in Q_{2M}\}$$
denote the first time a slow phase begins and where a node $X_i$ enters $Q_{2M}$. Provided that $\Ts_i < \ff$, we let $\Vs_i$ denote the displacement vector associated with this slow phase.
 We henceforth let
$$\Xsam := \{X_i \in \Xsaw\co \Ts_i < T \}$$
denote the $T$-self-avoiding nodes with a slow waypoint in $Q_{2M}$ and show that $\FF := \s(\Xsam)$ encodes the desired persistent information.

\subsubsection{Conditional decorrelation}
After conditioning on $\Xsam$, the covariance at two different time instants only pertains to the self-avoiding nodes crossing the $M$-box in the fast phase. Now, these nodes do not induce correlations, as they form an ephemeral background process of points that may be observed for a single instant, but are afterwards never seen again.
%
%
\bepr[Decay of correlations]
  \label{decorr_prop}
	  Let $s \ne t \le 1$. Then, 
		\begin{align}
			  \label{decorr_eq}
					\lim_{T \to\ff}\E\big[\Cov\big[\pia s, \pia t\,|\, \Xsam \big]\big] = 0.
		\end{align}
	\enpr

\subsubsection{Computation of conditional expectation}
	Taking into account Proposition \ref{decorr_prop}, in what follows, we need to determine the convergence of 
	$$\int_0^1 \E[\pia t\,|\, \Xsam] \de_t \d t = \int_0^1 \E[\pia t\,|\, \Xsam(tT)] \de_t \d t,$$
	as $\tff$. That means controlling the limit during the fast phases and during the slow phases.

We begin by simplifying the fast phase, where the fast nodes are essentially integrated out. More precisely, we replace them by a homogeneous Poisson point process with intensity $\laf(t)$.

	%
	%
	\bepr[Integrating out fast nodes]
	\label{fast_prop}
	Let $g\co\itf$ be continuous.	  Then, 
	$$\int_0^1 \E[\pia t\,|\, \Xsam(tT)] g(t)\d t - \int_0^1 \th(\Xsam(tT); \laf(t)) g(t)\d t$$
	tends to 0 in $L^1$.
	\enpr
	The final step consists in relating the trajectories from $\Xsam$ with the limiting birth-death process $\Xs$ in $Q_M$. Since the nodes $X_i \in \Xsam$ are $T$-self-avoiding, they visit $Q_M$ at most in the slow phase starting at $\Ts_i$ and possibly the fast phases immediately preceding and succeeding it. Hence, we associate with $\Xsam$ the birth-death process 
	$$\Xbam := \{(X_i(\Ts_i), \Ts_i, \Vs_i)\co X_i \in \Xsam\}.$$

	The final piece entering the proof of Proposition \ref{Thm_1m} is to show that $\Xbam$ indeed converges in distribution to the limiting process $\Xss$ in $Q_M$. 

	\bepr[Convergence of slow process]
	\label{slow_prop}
Under the scalings \eqref{scale_la_eq} and \eqref{scale_t_eq}, as $\tff$ it holds that 
	$$\Xbam \xrd \Xss.$$
	\enpr

	%
	%
	After laying out the general road-map, we now formally conclude the proof of Proposition \ref{Thm_1m}.
	\bep[Proof of Proposition \ref{Thm_1m}]
	Let $\e > 0$ and $g\co\itf$ be continuous.
	First,
	$$\tp_T(g) = \int_0^1 \pia t g(t)\d t,$$
	so that by the Chebyshev inequality conditioned on $\Xsam$,
	\begin{align*}
		\P\Big(\Big|\tp_T(g) - \int_0^1 \E[\pia t\,|\,\Xsam]g(t)\d t\Big| > \e\Big) &\le \frac1{\e^2}\E\Big[\Var\Big(\int_0^1 \pia t\,\Big|\, \Xsam\Big)g(t)\d t\Big]\\
		&= \frac1{\e^2}\int_0^1\int_0^1 \E\big[\Cov\big[\pia s, \pia t\,|\, \Xsam \big]\big]g(s)g(t)\d s\d t.
	\end{align*}
	By Proposition \ref{decorr_prop}, the expression in the last line tends to 0 as $\tff$. Hence, it suffices to prove the distributional convergence of 
	$\int_0^1\E\big[\pia t\,|\, \Xsam \big]g(t) \d t$. Proposition \ref{fast_prop} reduces this task further to the convergence of $\int_0^1 \th(\Xsam(tT); \laf(t))g(t)\d t$. Finally, since the map 
	$$\vp \mapsto \int_0^1 \th(\vp; \laf(t))g(t)\d t$$
	is continuous in the vague topology, the continuous mapping theorem together with Proposition \ref{slow_prop} conclude the proof.
	\enp

\subsection{Infrastructure-augmented model}
\label{out_inf_sec}
Next, we implement the convergence blueprint for the infrastructure-augmented model.

\subsubsection{Identification of persistent information}
\paragraph{Dense sinks}
 As the typical node sees a new set of relevant sinks essentially at all macroscopic times, no information persists in the limit and $\FF = \{\es, \Om\}$ is trivial. 

\paragraph{Sparse sinks}
Since the sinks are distributed sparsely, the typical node sees essentially the same set during the entire time horizon. Hence, the persistent information $\FF = \s(\vXi{}\, )$ consists of all sinks.

\paragraph{Critical setting}
In addition to the sinks, now also the path of the typical node becomes relevant. Hence, $\FF = \s(X_0(\cdot), \vXi{}\, )$.

\subsubsection{Conditional decorrelation}
The key observation is that at any two macroscopically differing times, the typical node sees different other nodes in the relevant vicinity, thereby guaranteeing near-perfect decorrelation.
%
%
To ease notation, we write
	$$\Ec t := \{\vXut(tT) \in \Xi^k(tT)\}$$
	  for the event of $k$-hop connection at time $tT$. Then, we establish
 the following time-decorrelation property.
\bepr[Decay of correlations]
	\label{lem-1b}
	Let $0 \le s < t\le 1$. Then, under the scaling \eqref{scale_la_eq} and \eqref{scale_t_eq}, 
\begin{align}
	\label{decorr_eq}
	\lim_{T\to\ff}\E\big[\Cov\big[\one\{\Ec s\}, \one\{\Ec t\}\ba \FF \big]\big] = 0;
\end{align}
holds for $\FF = \s(\vXut, \vXi{}\, )$ and $\FF = \s(\vXi{}\, )$. If additionally, $\a < d/2$, then \eqref{decorr_eq} also holds for $\FF = \{\es, \Om\}$.
\enpr

%
%
Before establishing Proposition \ref{lem-1b} in Section~\ref{sec_proofs}, we elucidate how to deduce Theorem \ref{thm_1} for $\a < d/2$.
\bep[Proof of Theorem~\ref{thm_1}, $\a < d/2$]
Let $\e > 0$, and $g\co\itf$ be continuous. Then,  by Chebychev's inequality, 
\begin{align*}
	\P\Big(\Big|\tktb(g) - \P(\Ec 0) \int_0^1g(t)\d t\Big|>\e\Big)&\le \frac1{\e^2}\Var[\tktb(g)]\\
	 &= \frac1{\e^2}\int_0^1\int_0^1\Cov[\one\{\Ec s\},\one\{\Ec t\}]g(s)g(t)\d s\d t,
\end{align*} 
where we could exchange integrals due to the boundedness of the integrand. Then, the result follows from dominated convergence and Proposition \ref{lem-1b}.
\enp

\subsubsection{Computation of conditional expectation}
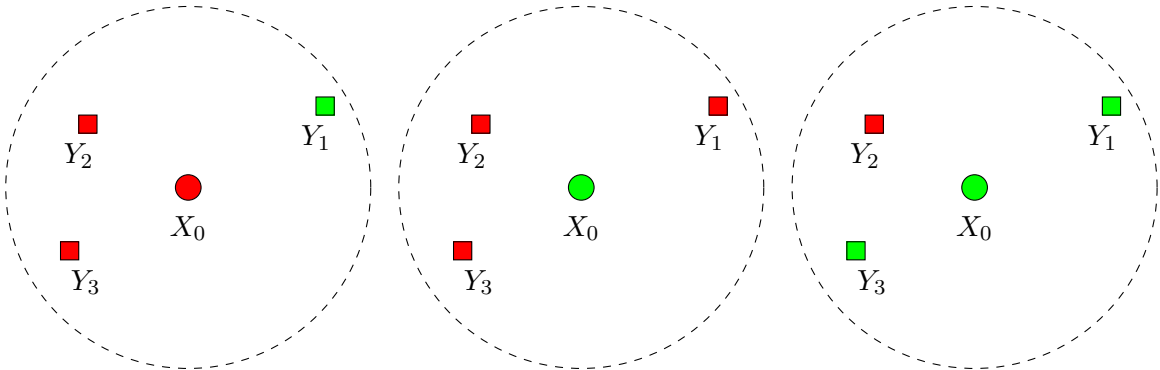
\begin{figure}[!htpb]
	\centering
	\begin{tikzpicture}[scale = 2.4]

	\draw[dashed] (0, 0) circle (1cm);

	\fill[red] (0,0) circle (2pt);
	\draw (0,0) circle (2pt);
	\coordinate[label=-90:{$X_{0}$}] (A) at (0,-.1);


	\fill[green] (.7, .4) rectangle (.8, .5);
	\fill[red] (-.6, .3) rectangle (-.5, .4);
	\fill[red] (-.7, -.4) rectangle (-.6, -.3);

	\draw (.7, .4) rectangle (.8, .5);
	\draw (-.6, .3) rectangle (-.5, .4);
	\draw (-.7, -.4) rectangle (-.6, -.3);

	\coordinate[label=-90:{$Y_{ 1}$}] (A) at (.7, .4) ;
	\coordinate[label=-90:{$Y_{ 2}$}] (A) at (-.6, .3) ;
	\coordinate[label=-70:{$Y_{ 3}$}] (A) at (-.7, -.4) ;
\end{tikzpicture}
\hspace{.1cm}
\begin{tikzpicture}[scale = 2.4]
	\draw[dashed] (0, 0) circle (1cm);

	\fill[green] (0,0) circle (2pt);
	\draw (0,0) circle (2pt);
	\coordinate[label=-90:{$X_{0}$}] (A) at (0,-.1);


	\fill[red] (.7, .4) rectangle (.8, .5);
	\fill[red] (-.6, .3) rectangle (-.5, .4);
	\fill[red] (-.7, -.4) rectangle (-.6, -.3);

	\draw (.7, .4) rectangle (.8, .5);
	\draw (-.6, .3) rectangle (-.5, .4);
	\draw (-.7, -.4) rectangle (-.6, -.3);

	\coordinate[label=-90:{$Y_{ 1}$}] (A) at (.7, .4) ;
	\coordinate[label=-90:{$Y_{ 2}$}] (A) at (-.6, .3) ;
	\coordinate[label=-70:{$Y_{ 3}$}] (A) at (-.7, -.4) ;

\end{tikzpicture}
\hspace{.1cm}
\begin{tikzpicture}[scale = 2.4]
	\draw[dashed] (0, 0) circle (1cm);

	\fill[green] (0,0) circle (2pt);
	\draw (0,0) circle (2pt);
	\coordinate[label=-90:{$X_{0}$}] (A) at (0,-.1);


	\fill[green] (.7, .4) rectangle (.8, .5);
	\fill[red] (-.6, .3) rectangle (-.5, .4);
	\fill[green] (-.7, -.4) rectangle (-.6, -.3);

	\draw (.7, .4) rectangle (.8, .5);
	\draw (-.6, .3) rectangle (-.5, .4);
	\draw (-.7, -.4) rectangle (-.6, -.3);

	\coordinate[label=-90:{$Y_{ 1}$}] (A) at (.7, .4) ;
	\coordinate[label=-90:{$Y_{ 2}$}] (A) at (-.6, .3) ;
	\coordinate[label=-70:{$Y_{ 3}$}] (A) at (-.7, -.4) ;

\end{tikzpicture}
\caption{Sparse sink case. The typical node and the sinks in the $k$-hop vicinity are represented by  disks and squares, respectively. Nodes in green are part of the giant component, nodes in red are not.}
	\label{shortFig}
\end{figure}

\paragraph{Sparse sinks}
In the sparse regime, the set of possible sinks that are within $k$-hop range does not change over the time horizon with a high probability. In particular, the randomness coming from this quantity remains in the limit. Figure \ref{shortFig} illustrates this phenomenon.

Hence, the probability that both, the typical node and at least one of the $N = \vXi{}(B_{k/\mu})$ relevant sinks, are in the unbounded connected component, equals $\th(\la)\big(1 - (1 - \th(\la))^N\big)$.
\paragraph{Critical setting}
 Since the movement of the typical node now is part of the persistent information, it must be taken into account. More precisely, the steps in the movement model are taken independently, so that by the invariance principle, in the limit, the typical node follows a path distributed according to a standard Brownian motion $\{W_t\}_{t \ge 0}$. Hence, at time $tT$, we obtain a connection with a conditional probability $\th(\la) \big(1- (1 - \th(\la))^{\vXi{}'(B_{1 / \mu}(W_t))}\big)$, where $\vXi{}'$ is a homogeneous Poisson point process with intensity $c_0$.


%
%

The major step in the proof is to show that conditioned on the sinks $\vXi{}$, the average $k$-hop connection time converges in distribution to the expression on the right-hand side in \eqref{short_eq}. 
For this, we show first that the percolation of the typical node and that of any sink can be considered independently. 
To make this precise, we let
	$$\Eii := \{\vXi i\in\CC(0)\text{ for some }\vXi i\in B_{k/\mu}\},$$
denote the event that some sink in $B_{k/\mu}$ is contained in the unbounded percolation cluster. 
In the critical case considered below, the movement over time matters, so that we also need to keep track of the location of the typical node. Therefore, we let
 $$\Ei t := \{\vXi{i}\in\CC(tT)\text{ for some }\vXi i\in B_{k/\mu}(\vXut(tT))\}$$
denote the event of seeing a sink in the unbounded component in the $k/\mu$-ball around $\vXut(tT)$.

%
%
\bepr[Asymptotic independence]
	\label{lem-2}
Let $g\co\itf$ be continuous and assume the scalings \eqref{scale_la_eq} and \eqref{scale_t_eq}. 
	\been
	\im Then,
		$$\E[\tktb(g)\ba \vXut, \vXi{}\, ] - \th(\la)\int_0^1 \P(\Ei t\ba\vXut, \vXi{}\, )g(t) \d t $$
		tends to 0 in $L^1$ as $\tff$.
	\im If $\a > d/2$, then
	$$\E[\tktb(g)\ba \vXi{}\, ] - \th(\la)\P(\Eii\ba\vXi{}\, ) \int_0^1g(t) \d t$$
		tends to 0 in $L^1$ as $\tff$.
		\enen
	\enpr
Finally, we establish the limit of the percolation probability $\P(\Eii\ba\vXi{}\, )$ for $\a > d/2$.

%
%
\bepr[Convergence of conditional $k$-hop connection probability]
	\label{lem-2b}
	Let $\a > d/2$. Then, under the scalings \eqref{scale_la_eq} and \eqref{scale_t_eq}, as $\tff$,
	$$\P(\Eii\ba\vXi{}\, ) \xrd 1 - (1 - \th)^N,$$
	where $N$ is a Poisson random variable with parameter $c_0|B_{1/\mu}|$. 
	\enpr

As soon as Propositions \ref{lem-2} and \ref{lem-2b} are available, the proof of Theorem \ref{thm_1} is completed through the second-moment method. 

%
%
\bep[Proof of Theorem \ref{thm_1}; $\a > d/2$]
Let $g\co\itf$ be continuous. With Propositions \ref{lem-2} and \ref{lem-2b} at our disposal, the task is to show that $\tktb(g)$ concentrates around the conditional mean, in the sense that $\tktb(g)- \E[\tktb(g)\ba \vXi{}\, ]$ tends to 0 in probability. 
	Arguing as in the proof of Theorem \ref{thm_2},  comparing $\tktb(g)$ and the random variable $\E[\tktb(g)\ba\vXi{}\, ]$ gives that
	\begin{align*}
		\P(|\tktb(g) - \E[\tktb(g)\ba\vXi{}\, ]|>\e)				&\le \frac1{\e^2}\int_0^1\int_0^1\E\big[\Cov[\one\{\Ec s\},\one\{\Ec t\}\,|\,\vXi{}\, ]\big]g(s)g(t)\d s\d t,
	\end{align*}
where we could exchange integrations due to the boundedness of the integrand. In particular, invoking dominated convergence and Proposition \ref{lem-1b} concludes the proof.
\enp

\begin{figure}[!htpb]
	\centering
	\begin{tikzpicture}[scale = 1.5]

	\draw[->, red, very thick] plot [smooth] coordinates {(0,0) (1,-.3) (2,.5) (3.5,-.3) (4,.2)} ;

	\draw[dashed] (0, 0) circle (1.5cm);
	\draw[dashed] (2,.5) circle (1.5cm);
	\draw[dashed] (4,.2) circle (1.5cm);

	\fill (0,0) circle (1pt);
	\fill (2,.5) circle (1pt);
	\fill (4,.2) circle (1pt);

	\coordinate[label=-90:{$X_{0}(0)$}] (A) at (0,0);
	\coordinate[label=-110:{$X_{0}(T/2)$}] (A) at (2.7,.9);
	\coordinate[label=-70:{$X_{0}(T)$}] (A) at (4,.2);

	\draw (0.7, 0.4) rectangle (0.8, 0.5);
	\draw (-0.6, 0.3) rectangle (-0.5, 0.4);
	\draw (-0.7, -0.4) rectangle (-0.6, -0.3);
	\draw (1.7, 0.5) rectangle (1.8, 0.6);
	\draw (2.7, -0.8) rectangle (2.8, -0.7);
	\draw (3.7, -0.4) rectangle (3.8, -0.3);

	\coordinate[label=-90:{$Y_{ 1}$}] (A) at (0.7, 0.4) ;
	\coordinate[label=-90:{$Y_{ 2}$}] (A) at (-0.6, 0.3) ;
	\coordinate[label=-70:{$Y_{ 3}$}] (A) at (-0.7, -0.4) ;
	\coordinate[label=-70:{$Y_{ 4}$}] (A) at (1.3, 0.7) ;
	\coordinate[label=-70:{$Y_{ 5}$}] (A) at (2.7, -0.4) ;
	\coordinate[label=-70:{$Y_{ 6}$}] (A) at (3.7, -0.4) ;

\end{tikzpicture}
	\caption{Critical scaling. Dashed circles illustrate the $k$-hop connection range of the typical node at different times.}
	\label{critFig}
\end{figure}
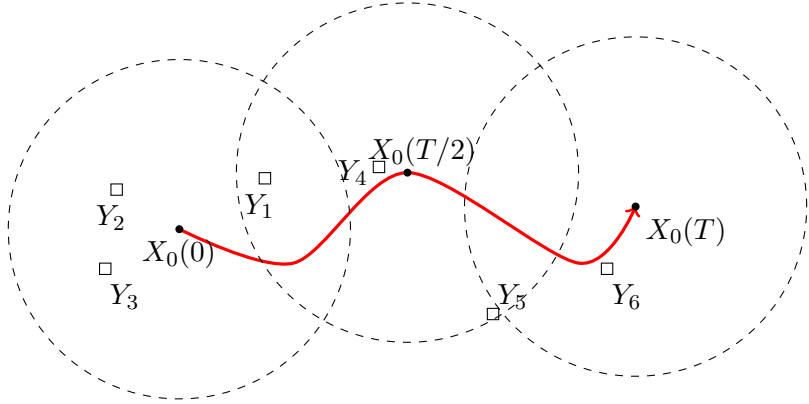

Finally, we study the critical scaling described in \eqref{crit_eq}. Recall that for sparse sinks the set of sinks that are within $k$-hop distance of the typical node remains the same during the entire time window. On the other hand, for dense sinks, we may assume that in practically every time point, we observe a fresh set of sinks. Loosely speaking, the critical regime interpolates between the two extremes. Figure~\ref{critFig} illustrates that as the typical node moves along its trajectory, it sees a finite number of new sinks, that often stay within $k$-hop reach for a substantial amount of time.

Due to the invariance principle, in the critical scaling, Brownian motion should appear in the limit. The key step to establish the description \eqref{crit_eq} is to identify the limit of the conditional probability.

%
%
\bepr[Convergence of conditional probabilities]
	\label{lem-2b-crit}
	Let $\a = d/2$ and $g\co\itf$ be continuous. Then, under the scalings \eqref{scale_la_eq} and \eqref{scale_t_eq}, as $\tff$,
	$$ \int_0^1g(t) \P(\Ei t\ba \vXut, \vXi{}\, ) \d t \xrd \int_0^1g(t)\big(1 - \big(1 - \th\big)^{\vXi{}'(B_{1 / \mu}(W_t))}\big) \d t,$$
where $\vXi{}'$ is a homogeneous Poisson process with intensity $c_0$. 
\enpr

%
%
Again, once Proposition \ref{lem-2b-crit} is established, the second-moment method enters the stage.
\bep[Proof of Theorem \ref{thm_1}; $\a = d/2$]
Propositions \ref{lem-2} and \ref{lem-2b-crit} reduce the claim to showing that for $s \ne t$ we have the correlation decay
$$\lim_{T\ua\ff}\E\big[\Cov[\one\{\Ec s\}, \one\{\Ec t\}\,|\,\vXut, \vXi{}\, ]\big] = 0,$$
	which holds by Proposition \ref{lem-1b}.
	\enp

\section{Proofs}\label{sec_proofs}
In Sections \ref{red_sec} and \ref{proof_inf_sec}, we present the proofs of the key technical auxiliary results for the two-scale mobility and the infrastructure-augmented model, respectively.
\subsection{Proofs for the two-scale mobility model}
\label{red_sec}
In Sections \ref{decorr_sec}--\ref{slow_sec}, we prove Propositions \ref{decorr_prop}--\ref{slow_prop}. In the remaining section, we abbreviate $\G=\G_0$ and do so accordingly for the waypoints and arrival times.

\subsubsection{Proof of Proposition \ref{decorr_prop}}
\label{decorr_sec}

%
%
As a preliminary step, we show that a typical trajectory is $T$-self-avoiding with high probability and  that it is highly unlikely for a node to visit a box at two different times during fast phases.

\bel[Non-recurrence of typical node]
\label{decorr_whp_lem}
Let $M \ge 1$ and $0 \le s < t \le 1$. Then,
\been
\im  $\lim_{\tff}\P(\G \text{ is $T$-self-avoiding}) = 1$, and
\im 
$\lim_{\tff}\P\big(\{|\G(sT) - \G(tT)|\le M\} \cap \{\{sT, tT\} \su \If\}\big) = 0.$
\enen
\enl

%
%
Now,  let
$$\Xsta := \{X_i \in X\co|\G_i(sT) - \G_i(tT)|\le M\text{ and }\{sT, tT\} \su \If\}$$
denote the nodes that are in their fast phase during times $sT$ and $tT$, and at those times, they are close to one another. Then, Lemma~\ref{decorr_whp_lem} together with the displacement theorem yields the following result for nodes visiting $Q_M$ at time $tT$.

	\bec[Non-recurrence of visiting nodes]
	\label{decorr_mtp_lem}
Let $M \ge 1$ and $0 \le s < t \le 1$. Then,
\been
\im 
	$\lim_{\tff}\E\big[X(tT)(Q_M) -  \Xsaw(tT)(Q_M)\big] = 0$, and
\im
	$\lim_{\tff}\E[\Xsta(tT)( Q_M)] = 0.$
\enen
	\enc

       %
        %
Now, we prove Proposition \ref{decorr_prop}. The key idea is to observe that thanks to Corollary~\ref{decorr_mtp_lem}, with high probability, fast nodes cannot be in $Q_M$ at two different times $s \ne t$. Thus, the independence property of the Poisson process establishes the desired mixing property. 
        \bep[Proof of Proposition \ref{decorr_prop}]
Let 
$$\Xmsta := \{X_i \in X\co X_i(tT) \in Q_M \text{ and }X_i(sT) \not \in Q_M\}$$ 
denote the family of all nodes that are in $Q_M$ at time $tT$ but not at time $sT$, and  set
$$\pimsta{} := \pi(\Xmsta(tT) \cup \Xsam(tT)).$$
        Then, by the independence property of the Poisson process,
        $\Cov\big[\pia s, \pimsta t\,|\, \Xsam \big] = 0,$
	so that it suffices to show that 
$\lim_{\tff} \E[\pia t - \pimsta t] = 0.$
By definition, the complement $X \sm \Xsam$ contains nodes whose trajectory is not $T$-self-avoiding and also nodes not entering $Q_M$ in a slow phase.  In particular, if $X_i \in X \sm \Xsam$ visits $Q_M$ both at times $s$ and $t$, then it is not $T$-self-avoiding or the visits must occur during the fast phase. In other words, $X_i(tT) \in (X(tT) \sm \Xsaw(tT)) \cap Q_M$ or $X_i(tT) \in\Xsta(tT) \cap Q_M$.
Hence, by the Markov inequality, it remains to show 
$$\lim_{\tff} \E\big[X(tT)(Q_M) - \Xsaw(tT)(Q_M)+ \Xsta(tT)(Q_M)  ] = 0,$$
so that applying Corollary~\ref{decorr_mtp_lem} concludes the proof.
        \enp

	%
	%
	\bep[Proof of Lemma \ref{decorr_whp_lem}] We split the proof into two parts and let $Q_M(x)$ denote the box of side length $M$, centered at $x\in\R^d$.
	
\noindent{\bf Part (1).}
Let $K > 0$ be arbitrary. Then, 
\begin{align*}\P(\G \text{ is not $T$-self-avoiding})\le&\sum_{\substack{j, j' \le K \\ j' \ne \{j, j + 1\}}}\P\big( Q_{4M}(\G(T_{2j})) \cap [\G(T_{2j' - 1}), \G(T_{2j'})] \ne \es\big)\\
& + \P(\#\{j \ge 0\co T_j^0\le 1\} > K).
\end{align*}
Since the second summand on the right-hand side becomes arbitrarily small for sufficiently large $K$, it suffices to show that for fixed $j, j'$ with $j' \not\in \{j, j + 1\}$ the first summand tends to 0 as $T \to \ff$. We prove the claim for $j' >  j + 1$, noting that the arguments for $j' < j$ are very similar. 

To that end, we take $L > 0$, set $U:= \G(T_{2j}) - \G(T_{2j ' -1 })$ and consider the decomposition
\begin{align*}
	\P\big( Q_{4M}(\G(T_{2j})) \cap [\G(T_{2j' - 1}), \G(T_{2j'})] \ne \es\big) \le  \P(|U| \le L) 
 + \P\big(|U| > L,\, Q_{4M}(U) \cap \R^+_0 V_{2j' - 1}  \ne \es \big).
 \end{align*}
We first show that the second summand becomes arbitrarily small for sufficiently large $L$. Indeed, suppose that there exists a point of the form $P = \G(T_{2j' - 1}) + t TV_{2j' - 1}$ with $P \in B_{4dM}(\G(T_{2j}))$ for some $t \ge 0$. Then, 
$$ |P - \G(T_{2j ' -1 })| \ge |U| - |P- \G(T_{2j})| \ge L - 4dM.$$ 
In particular, 
$$\f{V_{2j' - 1}}{|V_{2j' - 1}|}  = \f{P - \G(T_{2j ' -1 })}{|P - \G(T_{2j ' -1 })|}  \in B_{4dM / (L - 4dM)}\big(U\big).$$
Since we assumed $V_{2j' - 1}$ to be isotropic, ${V_{2j' - 1}}/{|V_{2j' - 1}|}$ is uniformly distributed on the unit sphere in $\R^d$. Therefore, the probability that the normalized increment ${V_{2j' - 1}}/{|V_{2j' - 1}|}$ is contained in $B_{4dM / (L - 4dM)}\big(U\big)$ tends to 0 as $L \to \ff$.

It remains to show that $\lim_{\tff}\P(|U| \le L) = 0$ for any fixed $L > 0$. To that end, we note that 
$$\G(T_{2j' - 1}) - \G(T_{2j})   = V_{2j} + TV_{2j + 1} + \cdots+ TV_{2j' -3}  + V_{2j' - 2}.$$
Then, by dominated convergence,
\begin{align*}
	\lim_{\tff}\P(|U| \le L) &= \E\big[\lim_{\tff} \one\big\{\big|V_{2j}/T + V_{2j + 1} + \cdots +V_{2j' -3}+ V_{2j' - 2}/T\big| \le L/T\big\}\big]\\
	&= \P(V_{2j + 1} + \cdots + V_{2j' - 3} = 0),
\end{align*}
which vanishes since we assumed $\k(\d v)$ to be absolutely continuous.

\noindent{\bf Part (2).} 
Similarly to Part (1), we fix $K > 0$ and then consider the decomposition
\begin{align*}
	\P\big(\{|\G(sT) - \G(tT)|\le M\} \cap \{\{sT, tT\} \su \If\}\big) &\le
	\P\big(\#\{j \ge 0\co T_j^0 \le 1\} > K\big)+ \sum_{ j < j'\le K} \P(E_{j, j', T})\\
	&\phantom \le+ \sum_{ j\le K} \P\big((t - s) T |V_{2j - 1}| \le M\big),
\end{align*}
where we set 
$$E_{j, j', T} := \Big\{\big|(\G(T_{2j' - 1}) + (tT - T_{2j' - 1})V_{2j' - 1}) - (\G(T_{2j - 1}) + (sT - T_{2j - 1})V_{2j - 1}) \big| \le M\Big\}.$$
As before, the first summand on the right-hand side becomes arbitrarily small for sufficiently large $K$.  Moreover, for fixed $K$, the final sum tends to 0 as $\tff$. Regarding the remaining term, we rely again on the limiting expression
$$\lim_{\tff}T^{-1}(\G(T_{2j' - 1}) - \G(T_{2j - 1})) = V_{2j - 1}  + V_{2j + 1} + \cdots + V_{2j' -3}.$$ 
Thus, 
we may again invoke dominated convergence  to deduce that
$$\lim_{\tff} \P(E_{j, j', T}) = \P\big((V_{2j - 1}  +  \cdots + V_{2j' -3}) + (t -T_{2j' - 1}^0) V_{2j' - 1} - (s -T_{2j - 1}^0) V_{2j - 1} = 0\big),$$
which vanishes due to the absolute continuity of $\k(\d v)$.

	\enp

\subsubsection{Proof of Proposition \ref{fast_prop}}
\label{fast_sec}
A key ingredient in the proof of Proposition \ref{fast_prop} is that at a given time $tT \le T$ most nodes in $\Xsam(tT)\cap Q_M$ are in their slow phase. To make this precise, we let
	$$\Xfa := \{X_i \in X\co t\in \If_i\}$$
	denote the family of nodes that are at time $tT$ in their fast phase. Then, this is the restriction to $Q_M$ of a homogeneous Poisson point process with intensity $\laf(t)$.

\bel[Nodes in $\Xsam\cap Q_M$ are slow]
\label{fast_mtp_lem}
Let $t \le 1$. Then,
$$\lim_{\tff}\E\big[(\Xfa(tT) \cap \Xsam(tT))(Q_M)\big] = 0.$$
\enl
Equipped with Lemma \ref{fast_mtp_lem}, we can now conclude the proof of Proposition \ref{fast_prop}.
        %
        %
        \bep[Proof of Proposition \ref{fast_prop}]
	First, for any continuous $g\co\itf$,
\begin{align*}
	&\E\Big[\Big|\int_0^1 \E[\pia t\ba \Xsam(tT)] g(t)\d t - \int_0^1 \th(\Xsam(tT); \laf(t)) g(t)\d t\Big|\Big]\\
	&\quad
	\le \int_0^1 \E\big[\big|\E[\pia t\ba \Xsam(tT)]  -  \th(\Xsam(tT); \laf(t))\big|\big] g(t)\d t,
\end{align*}
so that it suffices to prove $L^1$ convergence for fixed $t \le 1$. Now, if  nodes visit $Q_M$ in their slow phase and are not contained in $\Xsam$, then they necessarily fail to be self-avoiding. Thus, by the Markov inequality,
\begin{align*}	&\E\big[\big|\E[\pia t\ba \Xsam(tT)]  -  \E[\pi(\Xfa(tT) \cup \Xsam(tT))\ba \Xsam(tT)]\big|\big]\\
	&\quad\le  \E[\#\big((X(tT) \cap Q_M) \sm (\Xfa(tT) \cup \Xsam(tT))\big)] \\
	&\quad\le \E[X(tT)(Q_M) - \Xsaw(tT)(Q_M)].
\end{align*}
Hence, we can apply Corollary \ref{decorr_mtp_lem} and it remains to bound
\begin{equation}
	\begin{aligned}\label{y_cop_eq}
		&\E\big[\big|\E[\pi(\Xfa(tT) \cup \Xsam(tT))\ba \Xsam(tT)]  -  \th(\Xsam(tT); \laf(t))\big|\big] \\
		&\quad=\E\big[\big|\E[\pi(\Xfa(tT) \cup \Xsam(tT))\ba \Xsam(tT)]  -  \E[\pi(X' \cup \Xsam(tT))\ba \Xsam(tT)]\big|\big],
	\end{aligned}
\end{equation}
where we recall that $X'$ is a homogeneous Poisson point process with intensity $\laf(t)$ that is independent of $\Xsam(tT)$. To achieve this goal, it is critical to choose a wise representation for $X'$ on a suitably defined probability space. More precisely, by the independence property of the Poisson point process, the point process $\Xfa(tT) \sm \Xsam(tT)$ is a Poisson point process that is independent of $\Xsam(tT)$. However, since $\Xfa(tT)$ itself has the correct intensity measure, the intensity of $\Xfa(tT) \sm \Xsam(tT)$ is now to small. We can correct this by setting 
$$X': = (\Xfa(tT) \sm \Xsam(tT)) \cup X'',$$
where $X''$ is a point process that, when conditioned on $\Xsam(tT)$, is an independent copy of $\Xfa(tT) \cap \Xsam(tT)$.

Then, again the Markov inequality bounds the right-hand side of \eqref{y_cop_eq} by 
$$\E\big[\big(\Xfa(tT) \cap \Xsam(tT)\big)(Q_M)\big] +\E[X''(Q_M)]  = 2\E[\big(\Xfa(tT) \cap \Xsam(tT)\big) (Q_M)],$$
where the equality follows since we defined $X''$ to have the same distribution as $\Xfa(tT) \cap \Xsam(tT)$. Hence, an application of Lemma \ref{fast_mtp_lem} concludes the proof.
\enp

We end this section by proving Lemma \ref{fast_mtp_lem}.
\bep[Proof of Lemma \ref{fast_mtp_lem}]
By definition, nodes in $\Xsam$ visit $Q_M$ only in a single slow phase together with the fast phases immediately preceding and succeeding it. Hence, if $X_i$ is such that $X_i(tT)\in \Xsam(tT)\cap \Xfa(tT) \cap Q_M$, then $X_i$ changes phases in the time interval $[tT - M, tT + M]$. In other words, 
$$\E\big[\big(\Xfa(tT) \cap \Xsam(tT)\big)(Q_M)\big] \le \E\big[\#\{X_i(tT) \in Q_M:\, \{T_{i, j}\}_{j \ge0} \cap [tT - M, tT + M] \ne \es\}\big].$$
Now, by construction of the mobility model, we have $T_{j} = T T_{j}^0$, so that by the displacement theorem,
\begin{align*}
	\E\big[\big(\Xfa(tT) \cap \Xsam(tT)\big)(Q_M)\big] &\le \la |M|^d \P(\{T_j\}_{j \ge0} \cap [tT - M, tT + M] \ne \es)\\
	&= \la |M|^d \P(\{T_j^0\}_{j \ge0} \cap [t - M / T, t + M / T] \ne \es).
\end{align*}
But,
$$\limsup_{\tff}\P(\{T_j^0\}_{j \ge0} \cap [t - M / T, t + M / T] \ne \es)   \le \sum_{j \ge0}\P(T_j^0 = t) = 0,$$
since we assumed $\k(\d v)$ as absolutely continuous.
\enp

\subsubsection{Proof of Proposition \ref{slow_prop}}
\label{slow_sec}
The key step in the proof is to leverage the mass-transport principle \cite{mtp, mtp2} to determine the intensity measure $\mubd$. 
    \bel[Intensity of $\Xbam$]
        \label{slow_mtp_lem}
	The intensity measure $\mubd$ of $\Xbam$ is given by
	$$\mubd(\d x, \d t, \d v) = \la \one\{x \in Q_M\} \d x \E\Big[\sum_{j \ge 0}\de_{(T_{2j}^0, V_{2j})}(\d t, \d v) \one\{\G \text{ is $T$-self-avoiding}\}\Big].$$
        \enl
	\bep
We may restrict to sets of the form $Q \times B$, where $Q \su Q_M$ is a box and $B \su [0, 1] \times \R^d$ is Borel. Now, partition $\R^d$ into boxes $\{Q^k\}_{k \in K}$ of the same shape as $Q^0=Q$ and define the mass $\Phi(k, k')$ transported from $Q^k$ to $Q^{k'}$ as the 
number of nodes in $X_i \in \Xsaw$ starting in $Q^k$, having a slow waypoint in $Q^{k'}$ and conforming with condition $B$. More precisely, writing $T_i^{*,  k'}$ and $V_i^{*,  k'}$ for the arrival time at this waypoint and the associated displacement vector, 
$$\Phi(k, k') := \#\big\{X_i \in \Xsaw\co X_i(0) \in Q^k \text{ and } (T_i^{*,  k'}, V_i^{*,  k'}) \in B\big\}.$$
In particular,
$$\mubd(Q\times B) = \E\Big[\sum_{k \in K} \Phi(k, 0)\Big].$$
On the other hand, if $X_i \in \Xsaw$, then every $T_{i, 2j} \le T$ is of the form $T_i^{*,  k'}$ for a unique $k'$, so that 
$$\E\Big[\sum_{k' \in K} \Phi(0, k')\Big] = \la |Q| \E\Big[\sum_{j \ge 0}\one\{((T_{2j}^0, V_{2j}) \in B\} \one\{\G \text{ is $T$-self-avoiding}\}\Big].$$
Hence, applying the mass-transport principle concludes the proof.
      \enp
      Having computed the intensity $\mubd$, we now conclude the proof of Proposition \ref{slow_prop}.
%
%
\bep[Proof of Proposition \ref{slow_prop}]
Since both $\Xbam$ and $\Xss$ are Poisson point processes, it suffices to prove convergence of the intensity measure. Let $Q \su Q_M$ and $B \su [0, 1] \times \R^d$ Borel subsets. Then, by Lemma \ref{slow_mtp_lem}
	$$\mubd(Q \times B) =\la |Q| \E\big[\#\{j \ge 0\co (T_{2j}^0, V_{2j}) \in B\} \one\{\G \text{ is $T$-self-avoiding}\}\big].$$
By Lemma \ref{decorr_whp_lem}, the restriction of being $T$-self-avoiding disappears as $\tff$, so that indeed 
	$$\lim_{\tff}\mubd(Q \times B) =\la |Q| \E\big[\#\{j \ge 0\co (T_{2j}^0, V_{2j}) \in B\}\big] = \la|Q| (\nu \otimes \k)(B),$$
	as asserted.
\enp

\subsection{Proofs for the infrastructure-augmented model}
\label{proof_inf_sec}
In Sections \ref{dec_inf_sec}--\ref{crit_sec} we prove Propositions \ref{lem-1b}--\ref{lem-2b-crit}. We assume throughout that $k$, $\vli$ and $T$ are coupled according to \eqref{scale_la_eq} and \eqref{scale_t_eq}.
%
%

%
%
\subsubsection{Proof of Proposition \ref{lem-1b}}
\label{dec_inf_sec}
To prove decay of correlations, we rely on an implication of the shape theorem from \cite{bhopPerc}. This consequence reduces the $k$-hop connection event to finding a percolating sink in a $(k/\mu)$-ball. More precisely, for the convenience of the reader, we compress \cite[Lemmas 6 and 7]{bhopPerc} into a single result, where we let
\begin{align*}
	E^\ff_{t, T} := \{&\vXut(tT) \in \CC(tT)\text{ and }\vXi{ i} \in \CC(tT)\text{ for some }\vXi{ i}\in B_{k/\mu}(\vXut(tT))\}.
\end{align*}
be the event that at time $tT$ both $\vXut(tT)$ and some $\vXi{ i}\in B_{k/\mu}(\vXut(tT))$ are in the unbounded connected component $\CC(tT)$.

%
%
\bel[Consequence shape theorem; \cite{bhopPerc}]
\label{shape_lem}
Let $t \le 1$. Then,
$$	\lim_{T\to\ff}\P(\Ec t\De E^\ff_{t, T}) = 0.$$
\enl

Next, we approximate the true percolation events by percolation outside an $M$-box. That is, we put
\begin{align*}
	E_{t, T}^M := \{&\vXut(tT) \in \CCm(\vXu{}(tT))\text{ and }\vXi{ i} \in \CCm(\vXu{}(tT))\text{ for some }\vXi{ i}\in B_{k/\mu}(\vXut(tT))\}.
\end{align*}
where $\CCm(\vXu{}(tT))$ denotes the set of points $x \in \R^d$ that can leave their $Q_M$-neighborhood by relaying via nodes in $\vXu{}(tT)$. More precisely, there exist nodes $\vXu{i_1}(tT), \dots, \vXu{i_m}(tT)$ such that the following three conditions are satisfied, $|x - \vXu{i_1}(tT)|\le r$, $|\vXu{i_{j - 1}}(tT) - \vXu{i_j}(tT)| \le r$ for all $j \le m$, and $B_r(\vXu{i_m}(tT)) \cap  Q_M(x)^c\ne\es$.

%
%
\bel[$M$-box approximation]
\label{mapprox_lem}
It holds that
$$
	\lim_{M \to \ff}\limsup_{\tff}\P(E^M_{0, T} \sm E^\ff_{0, T}) = 0.
	$$
\enl
Before we prove Lemma~\ref{mapprox_lem}, let us show how it can be used to prove Proposition \ref{lem-1b}.
%
%
\bep[Proof of Proposition \ref{lem-1b}]
	By Lemmas \ref{shape_lem} and \ref{mapprox_lem}, we may replace $\Ec s$ by $E_{s, T}^M$ and $\Ec t$ by $E_{t, T}^M$. Since the movement is Markovian, we may take $s = 0$. Let us consider the different cases individually. \\[2ex]
{\bf $\boldsymbol{\FF = \s(\vXut, \vXi{}\, )}$.} Introduce the set of relevant sinks together with the typical node as 
$$\vXui(tT) := \{\vXut(tT)\}\cup \big(\vXi{}\cap B_{k/\mu}(\vXut(tT))\big).$$
Moreover, also introduce the $M$-neighborhoods,
$$\vWui(tT) := \vXui{}(tT) \oplus Q_M.$$
	Then, define $E_{t, T}^{M, *}$ just as $E_{0, T}^M$ except that for forming connections in $\CCm$, we only allow nodes $\vXu i$ not contained in $\vWui(tT)$ at time $tT$. 
	In particular, by the independence properties of the Poisson point process of nodes, 
$$\Cov[\one\{E_{t, T}^{M, *}\}, \one\{E_{t, T}^M\} \ba \vXut, \vXi{}\, ] = 0.$$
Furthermore, under the event $E_{0, T}^M \sm E_{t, T}^{M,*}$ there exists a relevant node at two times simultaneously,  i.e.,
	$$E_{0, T}^M \sm E_{t, T}^{M,*} \su F_k^M,$$
	where 
	$$F_k^M := \{\vXu{ i}(0)\in \vWui(0)\text{ and }\vXu{ i}(tT)\in \vWui(tT) \text{ for some $i \ge 1$} \}.$$
	Then, the probability of $F_k^M$ can be bounded via the Mecke formula \cite[Theorem 4.4]{poisBook}. To that end, we let $\G$ denote a random walk trajectory started at the origin. Then,
	\begin{align*}
		\P(F_k^M) 
		&\le \la\int_{Q_M}\E\Big[\#\{(z, z') \in \vXui(0) \times \vXui(tT)\co(\G(tT) +x + z) \in Q_M(z') \}\Big]  \d x\\
		&= \la\int_{Q_M}\E\Big[\sum_{z \in \vXui(0) }\sum_{z' \in \vXui(tT)}\P\big(\G(tT) \in Q_M(z' - z - x)\ba \vXut, \vXi{}\, \big)\Big]  \d x\\
		&\le \la c_0 \int_{Q_M}\int_{B_{1/ \mu}}\E\Big[\sum_{z \in \{\vXut(0), y\}}\sum_{ z'\in \{\vXut(tT), y\}}\P\big(\G(tT) \in Q_M(z' - z - x)\ba \vXut\big)\Big]\d y\d x \\
		&+\la c_0^2 \int_{Q_M}\int_{B_{1/ \mu}}\int_{B_{1/ \mu}}\E\Big[\sum_{z \in \{\vXut(0), y\}}\sum_{ z'\in \{\vXut(tT), y'\}}\P\big(\G(tT) \in Q_M(z' - z - x)\ba \vXut\big)\Big]\d y'\d y \d x .
	\end{align*}
	Now, by the central limit theorem, $\G(tT) / \sqrt T$ converges in distribution to a Gaussian random variable $Z$, so that for every $x' \in \R^d$, 
	$$\lim_{\tff}\P\big(\G(tT) \in Q_M(x')\big)  = \P(Z = o) = 0.$$
	Hence, an application of dominated convergence concludes the proof for ${\FF = \s(\vXut, \vXi{}\, )}$.\\[2ex]
{\bf $\boldsymbol{\FF = \s(\vXi{}\, )}$.} 
In this case, by the law of total covariance, 
\begin{align*}
	\E\big[\Cov\big[\one\{E_{0, T}^M\}, \one\{E_{t, T}^M\}| Y \big]\big] &=  \E\big[\Cov\big[\P\big(E_{0, T}^M|\vXut, \vXi{}\, \big), \P\big(E_{t, T}^M|\vXut, \vXi{}\, \big)|\vXi{}\, \big]\big]\\
	&\qquad + \E\big[\Cov\big[\one\{E_{0, T}^M\}, \one\{E_{t, T}^M\}| \vXut, \vXi{}\, \big]\big],
\end{align*}
where the first summand on the right-hand side vanishes       since $\P\big(E_{0, T}^M\ba\vXut, \vXi{}\, \big)$ is $\vXi{}$-measurable.
		Hence, invoking the result for $\FF = \s(\vXut, \vXi{}\, )$ concludes the proof.\\[2ex]
{\bf $\boldsymbol{\FF = \{\es, \Om\}}$.}  Note that the conditional probability $\P(E_{0, T}^M\ba \vXut)$ is almost surely constant. Hence, by the law of total covariance it suffices to prove that 
	\begin{align*}
	\lim_{\tff}\E\big[\Cov[\one\{E_{0, T}^M\}, \one\{E_{t, T}^M\} \ba \vXut]\big] = 0.
\end{align*}
Combining the result for $\FF = \s(\vXut, \vXi{}\, )$ with the law of total covariance, allows to further reduce the problem to proving that
	\begin{align*}
		\lim_{\tff}\E\big[\Cov[\P(E_{0, T}^M\ba \vXut, \vXi{}\, ), \P(E_{t, T}^M\ba \vXut, \vXi{}\, ) \ba \vXut]\big] = 0.
\end{align*}

Now, $\P(E_{0, T}^M\ba \vXut, \vXi{}\, )$ and $\P(E_{t, T}^M\ba \vXut, \vXi{}\, )$ depend only on $\vXi{} \cap B_{k/\mu}$ and $\vXi{} \cap B_{k/\mu}(\vXut(tT))$, respectively. Hence, under the event
$E_k := \{|\vXut(tT)|\ge 3k/\mu\}$
the independence of the Poisson point process of sinks yields that 
	\begin{align*}
		\lim_{\tff}\E\big[\Cov[\P(E_{0, T}^M\ba \vXut, \vXi{}\, ), \P(E_{t, T}^M\ba \vXut, \vXi{}\, ) \ba \vXut]\big] = 0.
\end{align*}

	To show that the probability of the event $E_k$ tends to 1 as $\tff$, we leverage again that by the central limit theorem the random vector $\vXu{ o}(tT)/\sqrt T$ converges in distribution to a centered Gaussian random variable $Z$. Therefore,
	\begin{align*}
		\lim_{\tff}\P(E_k^c) = \lim_{\tff}\P\big(|\vXu{ o}(tT)/\sqrt T| \le 3k^{1- d/(2\a)}/\mu\big) = \P(Z = o) = 0,
	\end{align*}
as asserted.
\enp
	%
	%
	As before, we write $\th^M$ for the probability to percolate beyond an $M$-box.
	\bep[Proof of Lemma \ref{mapprox_lem}]
	First, note that $E_{0, T}^M \sm E_{0, T}^\ff \su \ekmf$, where
	$$\ekmf:= \{o \in\CCm(\vXu{}(0)) \sm \CC(0)\} \cup \{\vXi i \in \CCm(\vXu{}(0)) \sm \CC(0)\text{ for some }\vXi i\in B_{k/\mu}(\vXut(tT))\}$$
	denotes the event that some relevant device reaches outside its $M$-neighborhood without being in the unbounded connected component. Then, by the Mecke formula \cite[Theorem 4.4]{poisBook}, 
	$$\P(\ekmf) \le (1 + \vli |B_{k/\mu}|) \P\big(o \in\CCm(\vXu{}(0)) \sm \CC(0)\big) = (1 + c_0|B_{1/\mu}|)(\th^M(\la) - \th(\la)).$$
	The right-hand side does not depend on $k$ and tends to 0 as $M \to \ff$.
	\enp

%
%
\subsubsection{Proof of Proposition \ref{lem-2}}
In the proofs of Propositions \ref{lem-2} and \ref{lem-2b} it is useful to consider the $M$-approximation
\begin{align*}
		\Emi := \{\vXi i \in \CCm  \text{ for some }\vXi i\in B_{k/\mu}\},
\end{align*}
where to ease notation, we write $\CCm$ instead of $\CCm(\vXu{}(0))$.

%
%
\bep[Proof of Proposition \ref{lem-2}]\phantom a\\[1ex]
{\bf Part (1).}	To prove the claim, note that by time-stationarity of the random-walk model,
	\begin{align*}
		\E\big[\big|&\E[\tktb(g)\ba \vXut, \vXi{}\, ] -\th(\la)\int_0^1 g(t)\P(\Ei t\ba\vXut, \vXi{}\, ) \d t\big|\big] 
		\\&= \E\big[\big|\int_0^1g(t)\P(\Ec t\ba \vXut, \vXi{}\, ) - \th(\la) g(t) \P(\Ei t\ba \vXut,\vXi{}\, )\d t \big| \big]\\
		&\le \int_0^1g(t)\E\big[\big|\P(\Ec t\ba \vXut, \vXi{}\, ) - \th(\la) \P(\Ei t\ba \vXut,\vXi{}\, )\big| \big]\d t\\
		&= \E\big[\big|\P(\Ec0\ba \vXi{}\, ) -\th(\la) \P\big(\Eii\ba\vXi{}\big)\big| \big]\int_0^1 g(t) \d t.
	\end{align*}
	Now, by Lemma \ref{mapprox_lem},
	$\lim_{M \ua \ff}\limsup_{T \ua \ff}\P(E_{0, T}^M \sm \Ec0) = 0$
	 and the event $\Ekm0$ decomposes as 
	$$\Ekm0 =  \{o \in \CCm\}\cap \Emi.$$
Since the nodes $\vXu{}$ form a Poisson point process, we conclude that under the event $\{\vXi{} \cap Q_{2M}= \es\}$
  the conditional probability $\P(\Ekm0\,|\,\vXi{}\, )$ factorizes as
	$$\P(\Ekm0\,|\,\vXi{}\, ) = \P(o \in \CCm)\P(\Emi\,|\,\vXi{}\, )= \th^M(\la)\P(\Emi\,|\,\vXi{}\, ).$$
	 Finally, we conclude the proof by noting that the event $\{\vXi{} \cap Q_{2M}= \es\}$ occurs with high probability, because
	 $$\P(\vXi{} \cap Q_{2M}= \es) = \exp(-\vli(2M)^d) = \exp(-T^{-\a}(2M)^d)$$ 
	 tends to 1 as $\tff$.
	 
\noindent	 
{\bf Part (2).}	Using the first part, it suffices to show that $\lim_{\tff} \P(\Eii \De \Ei t) = 0$ for every $t \le 1$. To that end, we leverage the Markov inequality to see that 
		\begin{align*}
			\P(\Eii \De \Ei t) \le \E\big[\vXi{} \big(\big(B_{k/\mu} \De B_{k/\mu}(\vXut(tT))\big)\big) \big] = c_0 \E\big[\big|B_{1/\mu} \De B_{1/\mu}(\vXut(tT)/k)\big|\big].
		\end{align*}
Since $\vXut(tT)k^{-1} = (\vXut(tT)/\sqrt T) k^{d/(2\a) - 1}$ and $\a > 2d$, the central limit theorem implies that the right-hand side in the above display tends to 0 as $T \to \ff$.
\enp
%
%
\subsubsection{Proof of Proposition \ref{lem-2b}}
\bep[Proof of Proposition \ref{lem-2b}]
	As in the proof of Proposition \ref{lem-1b}, we can relax the requirement of finding paths to infinity to finding paths leaving an $M$-box for some large $M > 0$. That is, we want to show that 
	$$ \P(\Emi\ba\vXi{}\, ) \xrd 1 - \big(1 - \th^M(\la)\big)^N.$$
Now, let 
$$G_k := \{|\vXi i - \vXi j| \ge k^{d/(2\a)} \text{ for all $\vXi i \ne \vXi j \in  B_{k/\mu}$}\}$$
denote the event that all sinks in $ B_{k/\mu}$ have Euclidean distance at least $k^{d/(2\a)}$.  Then, for sufficiently large $T$ and for $Y\in G_k$, 
\begin{align*}
	1 - \P(\Emi \,|\, \vXi{}\, )&=\prod_{\vXi{ i}\in B_{k/\mu}}\big(1 - \P(\vXi i \in \CCm\ba \vXi{}\, )\big)
	=\big(1 - \th^M(\la)\big)^{\vXi{} (B_{k/\mu})}.
\end{align*}
Noting that $\vXi{}(B_{k/\mu})$ is a Poisson random variable with parameter $\vli |B_{k/\mu}| = c_0 |B_{1/\mu}|$, it suffices to show that the events $G_k$ occur with high probability. To that end, we apply the Mecke formula to deduce that 
\begin{align*}
	\P(G_k^c) \le \vli^2 \int_{B_{k/\mu}}\int_{B_{k/\mu}}\one\{|y - y'| \le k^{d/(2\a)} \}\d y \d y'	= c_0^2 \int_{B_{1/\mu}}\big|B_{k^{d/(2\a) - 1}}(y)\big|\d y.
\end{align*}
Since $\a > d/2$, the latter expression tends to 0 as $\tff$.
\enp

%
%
\subsubsection{Proof of Proposition \ref{lem-2b-crit}}
\label{crit_sec}
\begin{proof}[Proof of Proposition \ref{lem-2b-crit}]
	Again, as in the proof of Proposition \ref{lem-2}, we can relax the requirement of finding paths to infinity to finding paths leaving an $M$-box for some large $M > 0$. That is, we claim that in distribution,
	$$ \int_0^1 g(t)\big(1 - \P(\Emmi t\, |\, \vXut, \vXi{}\, )\big) \d t \xrd \int_0^1g(t)\big(1 -\th^M(\la)  \big)^{Y' (B_{1 / \mu}(W_t))} \d t. $$
For this, let first $G_{T, M}$ denote the high-probability event that all sinks within distance $k/\mu$ of the path of $\vXut$ have distance at least $2M$ from one another.  Then, under $G_{T, M}$,
	$$1- \P(\Emmi t\ba \vXut, \vXi{}\, )   = (1 -\th^M(\la) )^{Y (B_{k / \mu}(\vXut(tT)))} = (1 -\th^M(\la) )^{Y' (B_{1 / \mu}(\vXut(tT) / k))},$$
	where $Y' = (Y/k)$ is now a homogeneous Poisson process with intensity $\vli k^d = c_0$. Taking the integral, we highlight the dependence on $Y'$ and $\{\vXut(tT)/k\}_{t \le 1}$ by writing 
	$$F(Y', \{\vXut(tT)/k\}_{t \le 1}) := \int_0^1 g(t)(1 -\th^M(\la) )^{Y' (B_{1 / \mu}(\vXut(tT) / k))}.$$
	Now, let $f$ be any bounded Lipschitz function of Lipschitz constant 1. Then, we want to show that 
	$$\lim_{T \ua \ff}\E[f(F(Y', \{\vXut(tT)/k\}_{t \le 1}))] = \E[f(F(Y', \{W_t\}_{t \le 1}))].$$
	To achieve this goal, by the invariance principle, it suffices to show that  the mapping
	$\g \mapsto \E[f(F(Y', \g))]$ is continuous in the Skorokhod topology outside a zero-set with respect to the distribution of Brownian motion. Hence, we fix $\e > 0$ and a continuous trajectory $\g$. First, since $\g$ is  continuous, it suffices to prove the claim with the Skorokhod norm replaced by the sup-norm. Now, let $\{\g_n\}_n$ be a sequence of right-continuous trajectories. Then, by the Lipschitz assumptions, 
	\begin{align*}\big|\E[f(F(Y', \g_n))] - \E[f(F(Y', \g))]\big| &\le \int_0^1 \E\big[|(1 -\th^M(\la) )^{Y' (B_{1 / \mu}(\g_n(t)))} - (1 -\th^M(\la) )^{Y'(B_{1 / \mu}(\g(t))) }|\big] \d t\\
		&\le \int_0^1 \E\big[Y'\big(B_{1 / \mu}(\g_n(t)) \Delta B_{1 / \mu}(\g(t))\big)\big] \d t\\
		&= c_0\int_0^1 \big|B_{1 / \mu}(\g_n(t)) \Delta B_{1 / \mu}(\g(t))\big| \d t.
	\end{align*}
	Now, we conclude the proof by noting that the right-hand side tends to 0 as $\g_n \to \g$ in the sup-norm. 
\end{proof}

	\section*{Acknowledgements}
This research was supported by Orange S.A., France grant CRE G09292, the German Research Foundation under Germany's Excellence Strategy MATH+: The Berlin Mathematics Research Center, EXC-2046/1 project ID: 390685689, and the Leibniz Association within the Leibniz Junior Research Group on Probabilistic Methods for Dynamic Communication Networks as part of the Leibniz Competition.

\phantomsection
\addcontentsline{toc}{section}{References}
\bibliography{refs}
\bibliographystyle{alpha}

\end{document}